\documentclass[11pt]{article}

\usepackage{xcolor}
\usepackage{xspace}
\usepackage{graphicx}
\usepackage{subfigure}
\usepackage{setspace}
\usepackage[utf8]{inputenc}
\usepackage{amsmath,amssymb,amsthm,mathrsfs,amsbsy}
\usepackage{lineno}
\usepackage{siunitx}
\usepackage{todonotes}
\usepackage[labelfont=bf,small]{caption}
\usepackage{algorithm}
\usepackage{algpseudocode}
\usepackage[capitalize]{cleveref}
\usepackage{pgfplots}
\usepackage{authblk}
\usepackage[numbers,square,sort&compress]{natbib}

\DeclareUnicodeCharacter{2212}{−}
\usepgfplotslibrary{groupplots,dateplot}
\usetikzlibrary{patterns,shapes.arrows}
\pgfplotsset{compat=newest}

\newcommand{\bm}[1]{\ensuremath{\mathbf{#1}}}
\newcommand{\bs}[1]{\ensuremath{\boldsymbol{#1}}}
\newcommand{\Vhat}{\hat{\bm V}}

\newtheorem{remark}{Remark}

\begin{document}

\title{CUR for Implicit Time Integration of Random Partial Differential Equations on Low-Rank Matrix Manifolds}

\author[1]{Mohammad Hossein Naderi}
\author[1]{Sara Akhavan}
\author[1]{Hessam Babaee\thanks{Corresponding author. Email:h.babaee@pitt.edu.}}
\affil[1]{\footnotesize Department of Mechanical Engineering and Materials Science, University of Pittsburgh}

\date{}
\maketitle

\begin{abstract}
Dynamical low-rank approximation allows for solving large-scale matrix differential equations (MDEs) with significantly fewer degrees of freedom and has been applied to a growing number of applications. However, most existing techniques rely on explicit time integration schemes. In this work, we introduce a cost-effective Newton's method for the implicit time integration of stiff, nonlinear MDEs on low-rank matrix manifolds. Our methodology is focused on MDEs resulting from the discretization of random partial differential equations (PDEs). Cost-effectiveness is achieved by solving the MDE at the minimum number of entries required for a rank-\(r\) approximation. We present a novel CUR low-rank approximation that requires solving the parametric PDE at \(r\) strategically selected parameters and \(\mathcal{O}(r)\) grid points using Newton's method. The selected random samples and grid points adaptively vary over time and are chosen using the discrete empirical interpolation method or similar techniques. The proposed methodology is developed for high-order implicit multistep and Runge-Kutta schemes and incorporates rank adaptivity, allowing for dynamic rank adjustment over time to control error. Several analytical and PDE examples, including the stochastic Burgers' and Gray-Scott equations, demonstrate the accuracy and efficiency of the presented methodology.
\end{abstract}

\section{\label{sec:Intro}Introduction}
Dynamical low-rank approximation (DLRA) of high-dimensional systems has recently  achieved remarkable success in the numerical simulation of a rapidly growing array of diverse engineering and scientific problems \cite{KL07}. DLRA was initially introduced in quantum chemistry for solving the Schrödinger equation, where it is known as the multiconfiguration time-dependent Hartree (MCTDH) method \cite{Beck:2000aa}.  Recently, DLRA and closely related variations of it have been employed to solve a growing list of diverse problems including stochastic partial differential equations (SPDEs) \cite{SL09,MN18,PB20}, Boltzmann transport and Vlasov equations \cite{EL19,BVT20,HW22,KS23}, turbulent combustion \cite{RNB21}, shallow water equations \cite{KKK23}, control \cite{BMS18}, detection of rare events \cite{FS16}, sensitivity analysis \cite{DCB22}, chemical kinetics \cite{NBGCL21,LB24,NLGBL24},  hydrodynamic stability analysis \cite{Babaee_PRSA,KBHH21,KNH24},  and deep learning \cite{SZK22,SZC23}.

 DLRA can be interpreted as reduced-order models (ROMs) for matrix or tensor differential equations (MDEs or TDEs), where the low-rank subspaces are expressed as time-dependent bases (TDBs) \cite{KL10}. DLRA offers two key advantages over ROMs based on static linear subspaces, which we refer to as static-ROMs. First, by evolving the subspaces, DLRA can adapt to instantaneous changes in the dynamics, overcoming some of the longstanding challenges of static-ROMs. This advantage is particularly important for problems characterized by slowly decaying Kolmogorov $n$-width, such as advection-dominated phenomena or turbulent dynamical systems. Second, building DLRA-based ROMs does not require the offline stage of collecting data to compute the low-rank subspace. Instead, evolution equations for TDBs are obtained via residual minimization directly from the full-order model (FOM) evolution equation. As a result, DLRA can be viewed as an \emph{on-the-fly} ROM that can adapt to changes of the dynamics.

There are many other dimension reduction techniques and ROMs based on TDB that differ from DLRA; see, for example, \cite{TSC18,P20,HD23,PR24}. These techniques are not discussed further in this paper, as the focus of the current work is on DLRA-based ROMs.

In the past decade, significant progress has been made in developing stable time integration schemes for the numerical solution of DLRA evolution equations, driven by the inherent instability of these equations in the presence of small singular values. This issue is particularly problematic because reducing low-rank approximation error requires increasing the rank, which involves resolving dimensions with smaller singular values. Consequently, encountering small singular values is more the rule than the exception, explaining the significant attention this issue has garnered.
To address this issue, a projector-splitting technique was introduced \cite{LO14}, which is based on splitting the projection onto the tangent space and is capable of handling small or zero singular values. The approach proposed in \cite{LO14} involves a backward time step, rendering it unstable for parabolic systems. In contrast, the robust basis update and Galerkin (BUG) integrators \cite{CL21,CKL22,CLS23} are stable in the presence of small or zero singular values, even for parabolic systems. The BUG integrators have first-order temporal accuracy. Recently, a second-order robust BUG integrator based on the midpoint rule was introduced \cite{CEKL24}. Higher-order time-integration schemes based on the rank-truncation of the time-discrete evolution equations have also been introduced \cite{KV19,RDV22,CL23,GBTT24}. All of these techniques are robust in the presence of small or zero singular values.

The aforementioned works have primarily focused on developing explicit time integration schemes, which are not suitable for stiff problems. Recently, a time integration scheme based on exponential time-differencing was presented, which mitigates the issue of stiffness in matrix differential equations where the source of stiffness is the linear operator \cite{CV23}.

Far fewer studies address fully implicit time integration schemes. The underlying difficulty lies in the computational cost constraints that any implicit time integration scheme must meet to be viable. In \cite{RV23}, an implicit algorithm for the time integration of tensor differential equations on low-rank tensor train manifolds is presented. The computational cost of the approach presented in \cite{RV23} for nonlinear tensor differential equations can be significant since the nonlinear map (e.g., non-polynomial nonlinearity) of a low-rank tensor can be a full-rank tensor. In \cite{NQE23}, an implicit time integration scheme for solving high-dimensional linear advection-diffusion partial differential equations on low-rank manifolds was presented.

What is currently lacking in the literature is the ability to use standard high-order implicit time integration schemes to solve arbitrarily nonlinear MDEs on low-rank matrix manifolds in a cost-effective manner. This paper aims to fill that gap. Our focus is on solving parametric PDEs on low-rank matrix manifolds, which are relevant to outer-loop applications such as optimization, inverse problems, and uncertainty quantification.

The key ingredients of the presented methodology are twofold: (i) a CUR low-rank approximation where the residual generated due to the low-rank approximation is set to zero at strategically selected rows and columns, and (ii) an efficient numerical method to solve the nonlinear equations at the selected entries using  Newton's method.  We use the above framework to develop standard implicit time integration schemes, including second-order Adams-Bashforth, Backward Differentiation Formula (BDF), and diagonally implicit Runge-Kutta (DIRK) of various orders. The methodology is applicable to nonlinear random PDEs.

The remainder of the paper is organized as follows: \cref{sec:Problem} introduces key mathematical preliminaries, including matrix differential equations and DLRA formulation. \cref{sec:Method} provides a detailed presentation of the proposed implicit TDB-CUR methodology. \cref{sec:DC} demonstrates the performance of the method on analytical test problems, stochastic Burgers' equation, and stochastic 2D Gray-Scott equations. Finally, \cref{sec:Conclusion} the concluding remarks are presented.


\section{\label{sec:Problem}Preliminaries}

\subsection{Notation}
We present the notation used in this paper. Vectors are denoted by lowercase bold font, e.g., $\bm v \in \mathbb{R}^n$ and matrices are denoted by uppercase bold font,  e.g., $\bm V \in \mathbb{R}^{n\times s}$  the low-rank matrix manifold and introduce the notations used throughout the paper. We use $\bm I_n$ to denote the identity matrix of size $n \times n$. Any low-rank matrix is shown with the hat symbol $( \hat{ \ \ } )$, e.g., $\hat{\bm V} \in \mathbb{R}^{n\times s}$.   

We use Matlab-style indexing where    $\bm{V}(\bm{p},:) \in \mathbb{R}^{r\times s}$ denotes a submatrix of  $\bm{V}  \in \mathbb{R}^{n \times s}$ where $\bm{p} = [p_1,p_2, \dots, p_r]$ contains the row indices, and similarly $\bm{V}(:,\bm{s})$ denotes the submatrix of $\bm V$ where $\bm{s}=[s_1,s_2,\dots, s_r]$ are the column indices. Finally, $\bm{A}^{\dagger} = (\bm{A}^{\mathrm T}\bm{A})^{-1}\bm{A}^{\mathrm T}$ denotes the Moore-Penrose pseudoinverse of matrix $\bm{A}$.

\subsection{Problem Setup}

We consider a nonlinear partial differential equation (PDE) with parametric randomness:
\begin{equation}\label{eq:FOM_Cont}
\frac{\partial v}{\partial t} = f(v;x,t,\xi),
\end{equation}
subject to appropriate initial and boundary conditions. In this equation, $v=v(x,t;\xi)$, where $v$  is a function of both spatial coordinates, denoted as $x$, and time denoted as  $t$, with an additional dependence on a set of random parameters $ \xi=\left(\xi_1, \xi_2 \cdots, \xi_d\right) \in \mathbb{R}^{d}$. The function $f(v;x,t,\xi)$ represents nonlinear spatial differential operators. We consider cases where $f(\cdot)$ is a nonlinear function of $v$, where the nonlinearity could be polynomial (e.g. quadratic, cubic, etc.) or non-polynomial (e.g. exponential, fractional, etc.).
Discretizing \cref{eq:FOM_Cont} using a method of lines results in the following nonlinear MDE:
\begin{equation}\label{eq:FOM1}
\frac{\mathrm d \bm V}{\mathrm d t} = \mathcal F(\bm V), \quad   t\in I=[0,T_f],
\end{equation}
where $I=[0, T_f]$ is the time interval and $\bm V: I \rightarrow \mathbb{R}^{n \times s}$ is a matrix. The function $\mathcal F(\bm V)$ is obtained by discretizing $f(v;x,t,\xi)$, with respect to both $x$ and $\xi$. The rows of $\bm V$ are associated with spatial degrees of freedom ($x$) and the columns of $\bm V$ are associated with different random samples ($\xi$). Therefore, each column of $\bm{V}$ corresponds to the solution of \cref{eq:FOM1} for a fixed choice of $\xi$ parameters.   We assume that boundary conditions have been incorporated in $\mathcal F(\bm V)$.

The columns of $\bm{V}$ can be determined independently, allowing for the separate computation of each column. In contrast, the rows of $\bm{V}$ are dependent, necessitating the simultaneous consideration of all rows, as the determination of one row requires knowledge of the values in the other rows. We consider sparse spatial discretization schemes, where each row requires the values of $p_a$ rows, with $p_a \ll n$. This is consistent with most discretization methods such as finite difference, finite volume, finite element, and spectral element methods, which naturally result in sparse row dependencies.

The computational cost of computing each column of the FOM scales with $\mathcal{O}(n^\beta)$, and solving the FOM for all $s$ columns scales as $\mathcal{O}(s n^\beta )$. Here, $\beta$ depends on the complexity of the solver used for the linear systems. In particular, $1 \leq \beta \leq 3$ where $\beta=3$ for direct solvers and $\beta=1$ for iterative solvers. 
 
\subsection{Dynamical Low-Rank Approximation (DLRA)}
DLRA provides an elegant mathematical framework for the time integration of MDEs on low-rank matrix manifolds \cite{KL07}. The key assumption is that  $\mathbf{V}(t)$ can be instantaneously approximated  well with a low-rank matrix $\hat{\bm V}(t)$ as shown below:
\begin{equation}\label{eqn:low-rank-approx}
\hat{\mathbf{V}}(t) = \bm{U}(t)\boldsymbol{\Sigma}(t)\bm{Y}(t)^{\mathrm T},
\end{equation}
 where $\mathbf{U} \in \mathbb{R}^{n\times r}$ and $\mathbf{Y} \in \mathbb{R}^{s\times r}$ are orthonormal spatial and parametric bases, $\boldsymbol \Sigma \in \mathbb{R}^{r\times r}$ is a full matrix, and $r \ll \min(n,s)$ is the rank. For brevity, we omit the explicit dependence on time. 
The idea behind TDB low-rank approximation is that the bases optimally evolve to minimize the residual due to the low-rank approximation. By substituting the low-rank approximation into the MDE, the residual is defined as:
\begin{equation}\label{eq:time-cont-res}
\mathbf{R}(t) = \frac{d}{dt}\left(\mathbf{U}\boldsymbol \Sigma \mathbf{Y}^{\top}\right) - \mathcal{F}\left(\mathbf{U}\boldsymbol \Sigma\mathbf{Y}^{\top}\right).  
\end{equation}
The evolution equations for $\bm U, \bm \Sigma $, and $\bm Y$ are obtained by minimizing the residual:
\begin{equation}
\mathcal{J}(\dot{\mathbf{U}}, \dot{\boldsymbol \Sigma}, \dot{\mathbf{Y}}) = \left\lVert \frac{d}{dt}\left(\mathbf{U}\boldsymbol \Sigma\mathbf{Y}^{\top}\right) - \mathcal{F}\left(\mathbf{U}\boldsymbol \Sigma\mathbf{Y}^{\top}\right) \right\rVert_{F}^{2},
\end{equation}
subject to orthonormality of $\mathbf{U}$ and $\mathbf{Y}$. The above constrained  minimization problem can be solved  using Riemannian optimization \cite{KL07} or by using Lagrange multipliers \cite{RNB21}. The resulting evolution equations are:
\begin{subequations}
\begin{align}
\dot{\boldsymbol \Sigma} &= \mathbf{U}^{\top}\mathbf{F}\mathbf{Y},
\label{eq:DBO_evol_S}\\
\dot{\mathbf{U}} &= \left(\mathbf{I}_n - \mathbf{U}\mathbf{U}^{\top}\right)\mathbf{F}\mathbf{Y}\boldsymbol \Sigma^{-1}, 
\label{eq:DBO_evol_U} \\
\dot{\mathbf{Y}} &= \left(\mathbf{I}_s - \mathbf{Y}\mathbf{Y}^{\top}\right)\mathbf{F}^{\top}\mathbf{U}\boldsymbol \Sigma^{-\top}, \label{eq:DBO_evol_Y}
\end{align}
\end{subequations}
where $\mathbf{F} = \mathcal{F}(\mathbf{U}\boldsymbol \Sigma\mathbf{Y}^{\top})$. This is equivalent to orthogonal projection of $\mathcal F(\bm U \boldsymbol \Sigma \bm Y^{\mathrm T} )$ onto the tangent space of  the low-rank matrix  manifold at $\hat{\bm V}=\bm U \boldsymbol \Sigma \bm Y^{\mathrm T}$. The DLRA framework has the potential to significantly reduce the computational cost of solving MDEs. 

Despite the potential computational cost savings of Eqs. \ref{eq:DBO_evol_S}-\ref{eq:DBO_evol_Y} for solving MDEs, several challenges remain for practical problems. In particular for nonlinear MDEs with non-polynomial nonlinearity, $\mathcal F(\hat{\bm V})$ is  full rank despite $\hat{\bm V}$ being low rank leading to computational costs that can surpass the cost of solving the FOM.  In such cases, the explicit computation of $\mathcal{F}(\mathbf{V})$ becomes necessary, incurring computational costs scaling as $\mathcal{O}(ns)$, similar to the FOM. Additionally, the efficient implementation of DLRA can be intrusive, requiring the substitution of the low-rank approximation into the FOM and the projection of $\mathcal{F}(t,\mathbf{V})$ onto the tangent space of the low-rank manifold. This necessitates careful implementation to avoid forming  matrices in the ambient space, i.e., matrices of size $n \times s$. Furthermore, \cref{eq:DBO_evol_U,eq:DBO_evol_Y} become unstable when $\boldsymbol{\Sigma}$ is singular or near-singular, which is problematic since retaining very small singular values is often necessary for an accurate approximation.

The stability issues due to small singular values have received significant attention in the past years and several successful methods have been proposed \cite{CL21,CKL22,CLS23,CEKL24,KV19,RDV22,CL23}. The issue of cost was also recently addressed \cite{NB23,DPNFB23}. All of these methods are based on explicit time integration. 

Implicit time integration of the DLRA evolution equations can be significantly more expensive than explicit time integration schemes. Even for linear MDEs, $\boldsymbol \Sigma$, $\bm U$, and $\bm Y$ are nonlinearly coupled.   In principle, it is possible to re-formulate DLRA evolution equations as a very large vector differential equation in the form of $\mathbf w =[\mbox{vec}(\boldsymbol \Sigma), \mbox{vec}(\bm U),\mbox{vec}(\bm Y)]^{^ \mathrm T}$, where $\bm w$ is a vector of size $\bm w \in \mathbb{R}^{N}$ with $N=r(n+s+r)$. The implicit time integration of vector $\bm w$ is computationally expensive because both $n$ and $s$ are typically very large. Moreover, similar to the explicit time integration, if $\mathcal F$ has non-polynomial nonlinearity, the matrix $\bm F=\mathcal F(\bm V)$ must be computed  and stored, which is computationally prohibitive for large $n$ and $s$.  

\section{\label{sec:Method}Methodology}
The DLRA can be viewed as minimizing the low-rank residual at the time-continuous level followed by the temporal discretization of the DLRA evolution equations. We refer to this approach as \emph{minimize-then-discretize}. We adopt an approach that reverses this order: we first perform the temporal discretization, then minimize the low-rank approximation residual at the time-discrete level. We refer to this approach as \emph{discretize-then-minimize}. Similar approaches have been use for solving MDEs and TDEs on low-rank manifolds \cite{KV19,RDV22}.  In this section, we present the time-discrete FOM and an implicit CUR-based residual minimization approach.

\subsection{\label{sec:FOM}Implicit Time Integration for the Full-Order Model}
First, we present implicit time integration for full-order nonlinear MDEs. To present the main ingredients of the methodology,  we first consider the implicit Euler method. The extension to implicit multistep and Runge-Kutta methods are presented later in this section. 

Let us consider the implicit Euler termpoal discretiation, where $\bm V^{k}$ is updated as follows:
\begin{equation}\label{eq:FOM}
\bm V^{k} = \bm V^{k-1} + \Delta t \mathcal{F}(\bm V^{k}),
\end{equation}
where $\Delta t$ is the time advancement. Since the target MDE is obtained by discretizing the random PDE given by Eq. \ref{eq:FOM_Cont}, when $\mathcal{F}$ is applied to $\mathbf{V}$, its action is understood to be column-wise: $\mathcal{F}(\bm{V})=[\mathcal{F}(\bm{v}_1) \  \mathcal{F}(\bm{v}_2) \  \dots \ \mathcal{F}(\bm{v}_{s}) ] $.
Newton's method is employed to solve the nonlinear MDE described by \cref{eq:FOM}, requiring iterative updates to the solution matrix as follows:
\begin{equation}\label{eq:Newton}
\bm V_{i+1}^{k} = \bm V_{i}^{k} + \delta \bm V_{i}^{k},
\end{equation}
where $\delta \bm V_{i}^{k} \in \mathbb{R}^{n \times s}$ is the Newton's correction matrix at the $i$-th iteration. By substituting \cref{eq:Newton} into \cref{eq:FOM}, we obtain:
\begin{equation}\label{eq:Newton2}
\bm V_{i}^{k} + \delta \bm V_{i}^{k} = \bm V^{k-1} + \Delta t \mathcal{F}(\bm V_{i}^{k} + \delta \bm V_{i}^{k}).
\end{equation}
Linearizing  the function $\mathcal{F}(\bm V_{i}^{k} + \delta \bm V_{i}^{k})$ around $\bm V_{i}^{k}$ results in:
\begin{equation}\label{eq:Taylor}
\mathcal{F}(\bm V_{i}^{k}(:,c) + \delta \bm V_{i}^{k}(:,c)) \approx \mathcal{F}(\bm V_{i}^{k}(:,c)) + \mathcal J(\bm V_{i}^{k}(:,c)) \delta \bm V_{i}^{k}(:,c), \quad c=1, 2, \dots, s,
\end{equation}
where $\mathcal J: \mathbb{R}^n \rightarrow \mathbb{R}^{n\times n}$ is the Jacobian function: $\mathcal J_{pq}(\bm v) = \partial [\mathcal F(\bm v)]_p/\partial  [\bm v]_q$, where $[\mathcal F(\bm v)]_p$ is the $p$-th entry of vector $\mathcal F(\bm v) \in \mathbb R^n$. It is important to note that each column of the matrix $\bm V_{i}^{k}$ has its own associated Jacobian matrix.
By substituting \cref{eq:Taylor} into \cref{eq:Newton2}, we obtain:
\begin{equation}\label{eq:Newton3}
\bm V_{i}^{k}(:,c) + \delta \bm V_{i}^{k}(:,c) = \bm V^{k-1}(:,c) + \Delta t \mathcal{F}(\bm V_{i}^{k}(:,c)) + \Delta t \mathcal J(\bm V_{i}^{k}(:,c)) \delta \bm V_{i}^{k}(:,c).
\end{equation}
By rearranging \cref{eq:Newton3}, we obtain:
\begin{equation}\label{eq:Newton4}
(\bm I -  \Delta t \mathcal J(\bm V_{i}^{k}(:,c))) \delta \bm V_{i}^{k}(:,c) = \bm V^{k-1}(:,c) - \bm V_{i}^{k}(:,c) + \Delta t \mathcal{F}(\bm V_{i}^{k}(:,c)).
\end{equation}
Now we can drop the superscript $k$ here and write \cref{eq:Newton4} as:
\begin{equation}\label{eq:Newton5}
\bm A_{i}^{(c)} \delta \bm V_{i}(:,c) =\bm b^{(c)}_i, \quad c=1, 2, \dots, s,
\end{equation}
where
\begin{align}
\bm A_i^{(c)}  &= \bm I - \Delta t \mathcal{J}(\bm V^k_i(:,c)), \\
\bm b^{(c)}_i &=\bm V^{k-1}(:,c) - \bm V_{i}^{k}(:,c) + \Delta t \mathcal{F}(\bm V_{i}^{k}(:,c)).
\end{align}
Here $\bm A_i^{(c)} \in \mathbb{R}^{n\times n}$ and $\bm b^{(c)}_i \in \mathbb{R}^n$. \cref{eq:Newton5} is a linear system of equations that is solved to find $\delta \bm V_{i}$ at each iteration and for each column of $\bm V_{i}$. The computational cost of solving the above linear system of equations for each sample using direct solvers scales as $\mathcal{O}(n^3)$. Since each sample can be solved independently, the total computational cost of solving the FOM scales as $\mathcal{O}(sn^3)$. For large sparse systems, iterative methods are used, resulting in a computational cost that scales as $\mathcal{O}(sn)$.

\subsection{Optimal Low-Rank Approximation with SVD}
In this section, we present an optimal low-rank approximation to compute the solution of the time-discrete FOM, given by  \cref{eq:FOM}, on a low-rank matrix manifold.  To this end, assuming that the solution at the previous time step is in the low-rank form, the solution at the next time-step can be obtained by using the FOM as shown below:
\begin{equation}\label{eq:FOM_TDB}
\bm V^{k} = \hat{\bm V}^{k-1} + \Delta t \mathcal{F}(\bm V^{k}).
\end{equation}
Despite utilizing the rank-$r$ approximation $\hat{\bm V}^{k-1}$ in \cref{eq:FOM_TDB}, the resulting matrix $\bm V^{k}$ is not of rank-$r$. To address this issue, we require a rank truncation operation to retract the solution back onto the rank-$r$ matrix manifold. A brute-force yet optimal method for rank truncation would be to perform SVD on $\bm{V}^k$:
\begin{equation}\label{eq:TD_SVD}
\hat{\bm V}_{best}^k = \texttt{SVD}(\bm V^k),
\end{equation}
where $\texttt{SVD}(\bm V^k) = \bm U_{best}^k \boldsymbol \Sigma_{best}^k {\bm Y^{k^{\mathrm T}}_{best}}$ denotes the rank-$r$ truncated SVD of the matrix $\bm V^k$. This approach provides an optimal low-rank approximation, and it does not involve the inversion of $\boldsymbol \Sigma$. Therefore, the time integration remains stable even in the presence of zero singular values.   However, this approach  does not offer a computational advantage for most problems. The challenge arises because the exact rank of $\bm{V}^k$ can be very large or $\bm{V}^k$ could be even full rank. Consider, for instance, the case where $\mathcal{F}(\cdot)$ is quadratically nonlinear, such as in the Burgers equation. If we start with a rank-$r$ matrix $\hat{\bm{V}}_1^k$ as the initial guess for Newton's iterations, the exact rank of $\mathcal{F}(\hat{\bm{V}}_i^k)$ will be $\mathcal{O}(r^2)$. Consequently, the rank of the matrix $\bm{B}_1 = [\bm{b}_1^{(1)}, \bm{b}_1^{(2)}, \dots, \bm{b}_1^{(s)}]$ is $\mathcal{O}(r^2)$. This implies that the rank of Newton's correction matrix $\delta \bm{V}_1$ in the first iteration is at least $\mathcal{O}(r^2)$. As a result, the rank of $\hat{\bm{V}}_2^k$ obtained from \cref{eq:Newton} is also $\mathcal{O}(r^2)$. The exact rank of $\hat{\bm{V}}_i^k$ increases exponentially with subsequent Newton iterations, meaning the rank of $\hat{\bm{V}}_{i+1}^k$ is $\mathcal{O}(r^{2i})$. If $\mathcal{F}(\cdot)$ has non-polynomial nonlinearity, then the Newton correction matrix is full rank in the first iteration. 

To prevent the rank growth during Newton's iteration, \cite{RV23} introduces a rank truncation after each iteration. However, as previously mentioned, for nonlinear MDEs with non-polynomial nonlinearity, the solution matrix becomes full rank after the first Newton iteration.

\subsection{Near-Optimal Low-Rank Approximation with CUR}\label{sec:sparse-sampling}

This section outlines an efficient approximation approach to overcome the computational challenges mentioned above. We present an implicit time integration algorithm based on a CUR low-rank approximation, which results in a collocation method to solve the residual minimization problem.

A CUR decomposition provides a rank-$r$ approximation of a matrix ${\bm{V}} \approx \bm{C} \tilde{\bm{U}} \bm{R}$, where $\bm{C} \in \mathbb{R}^{n \times r}$ and $\bm{R} \in \mathbb{R}^{r \times s}$ are composed of actual independent columns and rows from $\bm{V} \in \mathbb{R}^{n \times s}$. The matrix $\tilde{\bm{U}} \in \mathbb{R}^{r \times r}$ is selected to ensure that the CUR decomposition approximates ${\bm{V}}$ accurately. The choice of columns, rows, and the matrix ${\tilde{\bm{U}}}$ leads to different CUR decompositions, all of which results in different rank-$r$ matrices. In other words, the result of $\texttt{CUR}(\bm{V})$ is a low-rank matrix $\hat{\bm{V}}$. 

CUR decompositions, also known as pseudoskeleton or cross approximations, were first introduced in \cite{GTZ97}. CUR decompositions have been widely used in data analysis  \cite{MD09}. The popularity of CUR decompositions stems largely from the fact that they retain the actual columns and rows of the target matrix. This results in more interpretable low-rank approximations and better compression ratios for sparse matrices compared to SVD. Since $\bm{C}$ and $\bm{R}$ consist of actual data entries, they are fully interpretable and preserve the sparsity of the target matrix. In contrast, SVD generates singular vectors that are linear combinations of the columns and rows of the target matrix, often losing sparsity and becoming more difficult to interpret.  Cross low-rank approximations have also been extended to tensor low-rank approximations \cite{OT10,ACC21}. 

The accuracy of CUR low-rank approximation depends on the choices of columns and rows and how matrix $\tilde{\bm{U}}$ is computed. One useful CUR decomposition for computational purposes is:
\begin{equation}\label{eq:CUR}
    \texttt{CUR}(\bm V) = \bm V(:,\bm s) \bm V(\bm p, \bm s)^{-1}  \bm V(\bm p,:),
\end{equation}
where $\bm s =[s_1, s_2, \dots, s_r]$ and $\bm p =[p_1, p_2, \dots, p_r]$ are the indices of columns and rows, respectively. It is shown in \cite{GTZ97} that the accuracy of the CUR decomposition is related to the matrix volume, specifically the determinant of the intersection matrix formed by the selected rows and columns. Given that this selection problem is NP-hard, several heuristic algorithms have been proposed, including Maxvol \cite{GTZ97, GT01}, Cross2D \cite{T00}, leverage score \cite{MD09}, and the discrete empirical interpolation method (DEIM) \cite{SE16}.

Another factor that determines the accuracy of the CUR decomposition is how matrix $\tilde{\bm{U}}$ is computed. The DEIM CUR algorithm analyzed in \cite{SE16} is based on:   $\hat{\bm V} = \bm V(:,\bm s) \tilde{\bm{U}}  \bm V(\bm p,:)$, where $\tilde{\bm{U}}$ is obtained via the orthogonal projection of the target matrix onto the space spanned by the selected columns and rows:  $\tilde{\bm{U}} = (\bm C^{\mathrm T} \bm C)^{-1} \bm C^{\mathrm T}\bm V \bm R^{\mathrm T}(\bm R \bm R^{\mathrm T})^{-1}$, where $\bm C = \bm V(:,\bm s)$ and $\bm R = \bm V(\bm p,:)$. This CUR is more accurate than the CUR presented in \cref{eq:CUR}. However, the CUR based on orthogonal projection requires access to all entries of $\bm V$ as opposed to the CUR given in \cref{eq:CUR}, where $\tilde{\bm{U}} = \bm V(\bm p, \bm s)^{-1}$ and its computation requires access to only $r$ columns and rows of $\bm V$.

In \cite{DPNFB23}, a CUR algorithm is analyzed based on \cref{eq:CUR}, where $\bm{p}$ and $\bm{s}$ are determined by applying the DEIM algorithm to the exact or approximate left and right singular vectors of $\bm{V}$, respectively. The CUR DEIM low-rank approximation error is shown in \cite[Theorem 2.8]{DPNFB23} to be bounded by:
\begin{equation}\label{eq:CUR_bnd}
    \|\bm V - \hat{\bm V}\|_2 \leq c \hat{\sigma}_{r+1},
\end{equation}
where \(\hat{\sigma}_{r+1}\) represents the orthogonal projection error, defined as \(\hat{\sigma}_{r+1} = \max\{ \|(\bm I - \bm U \bm U^{\mathrm T}) \bm V \|_2, \|\bm V(\bm I - \bm Y \bm Y^{\mathrm T})\|_2\}\), and \(\|\cdot\|_2\) denotes the spectral norm. Here, $\bm{U} \in \mathbb{R}^{n \times r}$ and $\bm{Y} \in \mathbb{R}^{s \times r}$ are sets of orthonormal vectors, meaning $\bm{U}^{\mathrm{T}} \bm{U} = \bm{I}$ and $\bm{Y}^{\mathrm{T}} \bm{Y} = \bm{I}$, that span the same subspaces as $\bm{V}(:, \bm{s})$ and $\bm{V}(\bm{p}, :)^{\mathrm{T}}$, respectively.

When \(\bm U\) and \(\bm Y\) correspond to the exact singular vectors, \(\hat{\sigma}_{r+1}\) equals the \((r+1)\)-th singular value of \(\bm V\), representing the minimal error for a rank-\(r\) approximation. The amplification factor \(c \geq 1\) is dependent on the condition numbers of two matrices and is given by \(c = \min\{ \eta_r(1+\eta_c), \eta_c(1+\eta_r)\}\), where \(\eta_r = \|\bm U^{-1}(\bm p,:)\|_2\) and \(\eta_c = \|\bm Y^{-1}(\bm s,:)\|_2\).

According to \cite[Lemma 3.2]{CS10} \(\eta_r\) (or \(\eta_c\)) is bounded by:
\begin{equation}
    \eta_r \leq (1 + \sqrt{2n_1})^{r-1} \|\bm u_1\|_{\infty}^{-1},
\end{equation}
where \(\bm u_1\) is the first left singular vector, i.e., the first column of \(\bm U\). A similar bound applies to \(\eta_c\). This result demonstrates that \(\eta_r\) (or \(\eta_c\)) remains bounded regardless of the singular values, ensuring that the DEIM remains well-conditioned as \(r\) increases. Although the above error bound appears pessimistic, in practice, \(\eta_r\) and \(\eta_c\) are typically small. DEIM, in fact, selects the interpolation points using a greedy algorithm that minimizes \(\eta_r\) and \(\eta_c\).

The computation of CUR low-rank approximation using \cref{eq:CUR} is numerically unstable \cite{OT10}. A stable algorithm can be obtained by first performing QR decomposition on the selected columns:
\begin{equation}
\bm Q \bm R = \texttt{qr}(\bm V(:,\bm s)),
\end{equation}
where $\bm Q \in \mathbb{R}^{n \times r}$ is the orthonormal matrix and $\bm R^{r \times r}$. The low-rank matrix can be computed by interpolating every column $\bm V$ onto $\bm Q$ at DEIM-selected rows:
\begin{equation}
\hat{\bm V}  = \bm Q \bm Z,
\end{equation}
where $\bm Z = \bm Q(\bm p,:)^{\dagger} \bm V(\bm p,:)$. Here, we use the pseudoinverse of $ \bm Q(\bm p,:)$ for cases where the size of the selected rows is larger than $r$. This happens if row oversampling is used. Oversampling the rows improves the condition number of the DEIM algorithm by reducing $\eta_r = \|\bm U(\bm p,:)^{\dagger} \|$ and results in tighter error bounds in \cref{eq:CUR_bnd}. 

The expression of $\hat{\bm V}  = \bm Q \bm Z$ can be converted to the SVD of $\hat{\bm V} = \bm U \bs \Sigma \bm Y^{\mathrm T}$ as shown in \cite{DPNFB23}. The stable CUR algorithm is presented in \cref{alg:SCUR}.

\begin{remark}
Let $\bm V \in \mathbb{R}^{n\times s}$ be a matrix with a rank larger than $r$ and let $\hat{\bm V}$ be the CUR rank-$r$ approximation of $\bm V$ constructed according to \cref{alg:SCUR}. If row oversampling is not used, i.e., the number of selected rows and columns are equal, then $\hat{\bm V}(\bm p,:) = \bm V(\bm p,:)$ and $\hat{\bm V}(:,\bm s) = \bm V(:,\bm s)$. However, if row oversampling is used, i.e., the number of selected rows is larger than columns, then $\hat{\bm V}(\bm p,:) \neq \bm V(\bm p,:)$ and $\hat{\bm V}(:,\bm s) = \bm V(:,\bm s)$. 
\end{remark}

CUR decompositions are also particularly appealing for solving MDEs on low-rank matrix manifolds because they enable the development of accurate and stable algorithms that require computing only $r$ columns and rows of the matrix $\bm{V}$ to construct a low-rank approximation, rather than computing the entire matrix $\bm{V}$ as required for SVD. Recently, we developed CUR algorithms for cost-effective time integration of nonlinear MDEs on low-rank matrix manifolds \cite{NB23,DPNFB23}. Similar cross algorithms have also been developed for solving nonlinear TDEs on low-rank Tucker tensor and tensor train manifolds \cite{GKM18,GB24,D24,GBTT24}. These developments are limited to explicit time integration.

\subsubsection{Residual Collocation with CUR}
Replacing $\bm V^k$ with a rank-$r$ matrix $\Vhat^k$ in \cref{eq:FOM_TDB} generates a residual due to the low-rank approximation as follows:
\begin{equation}\label{eq:Residual}
\mathcal R(\hat{\bm V}^{k}) = \hat{\bm V}^{k} - \hat{\bm V}^{k-1} - \Delta t \mathcal{F}(\hat{\bm V}^{k}),
\end{equation}
where $\bm R =\mathcal R(\hat{\bm V}^{k}) \in \mathbb{R}^{n\times s}$ denotes the low-rank approximation residual matrix and in general it is a nonlinear map of $\hat{\bm V}^{k}$. This residual is the time-discrete counterpart of the time-continuous residual given by Eq. \ref{eq:time-cont-res}. The time-discrete variational principle can be posed as finding the optimal rank-$r$ approximation $\hat{\bm V}^k$ that minimizes the Frobenius norm of the residual, i.e., $\|\mathcal R(\hat{\bm V}^{k}) \|_F$.

In this work, we present an implicit CUR algorithm as follows:
\begin{equation}\label{eq:TD_CUR}
\hat{\bm V}^k = \texttt{CUR}(\bm V^k).
\end{equation}

We can interpret \cref{eq:TD_CUR} from the perspective of residual minimization. The CUR approach is equivalent to a residual collocation method where the residual is set to zero at $r$ strategically chosen rows and columns, i.e., the collocation entires. Specifically, this implies that $\bm{R}(\bm{p}, :) = \bm{0}$ and $\bm{R}(:, \bm{s}) = \bm{0}$, where $\bm{p}$ and $\bm{s}$ are vectors representing the row and column indices at which the residual is zero. In essence, this approach solves the full-order model (FOM) at selected rows and columns and it does not use projection to the tangent space as in DLRA.

 We use either DEIM \cite{CS10} or its variant QDEIM \cite{DG16} to find $\bm p$ and $\bm s$. As demonstrated in \cite{DG16}, the performance of the DEIM and QDEIM algorithms is very similar. For simplicity, we refer to the sampling algorithm as \texttt{DEIM}, with the understanding that both DEIM and QDEIM may be used.

A key advantage of using CUR in comparison to SVD is that it only requires computing $\bm{V}^k(\bm{p},:)$ (selected rows) and $\bm{V}^k(:,\bm{s})$ (selected columns) at each time step, avoiding full model evaluations. Computing $\bm{V}^k(\bm{p},:)$ and $\bm{V}^k(:,\bm{s})$ requires solving nonlinear systems, unlike the explicit CUR algorithm, where the solution at the selected rows and columns is computed directly from the solution at previous time steps. Another advantage of this approach is that it is stable in the presence of small or zero singular values.

The details of computing the sampled row and column information in the TDB-CUR framework are explained below.

\subsubsection{Computing the Columns}
To compute the selected columns $\mathbf{V}^{k}(:,\mathbf{s})$, we leverage the fact that for MDEs arising from the discretization of parametric PDEs, each column can be solved independently. To this end,  we compute $\mathbf{V}^{k}(:,\mathbf{s})$, which requires an independent implicit solve for each column. The column indices are obtained by applying the DEIM algorithm \cite[Algorithm 1]{CS10} to $\bm{Y}^{k-1}$, the matrix of right singular vectors of $\hat{\bm{V}}^{k-1} = \bm{U}^{k-1} \bs{\Sigma}^{k-1} \bm{Y}^{k-1^{\mathrm{T}}}$, which is already computed in the previous time step. Therefore, the selected columns vary at each time step according to:
\begin{equation}\label{eq:DEIM_s}
    \bm s = \texttt{DEIM}(\bm Y^{k-1}).
\end{equation}

We use Newton's method for the implicit time advancement of the selected columns, iteratively solving \cref{eq:Newton5} for $c = s_1, s_2, \dots, s_r$ until convergence is achieved. The proposed approach offers two computational advantages: (i) existing deterministic codes can be utilized in this step, as solving for each column requires time advancement of a deterministic code in a non-intrusive black-box fashion for a specific choice of parameter $\xi$; (ii) since the computation of the columns can be performed independently, this step is parallelizable.

\subsubsection{Computing the Rows}
To compute  $\mathbf{V}^k(\mathbf{p},:)$, we need to solve the following system of equations:
\begin{equation} \label{eq:LSoE}
    \bm A_{i}^{(c)}\delta \bm V_{i}(:,c)=\bm b_{i}^{(c)}
\end{equation}
at only $\bm p$ indices. 
The selected rows $\delta \bm V_{i}(\mathbf{p},c)$, can be computed as:
\begin{equation}
    \delta \bm V_{i}(\mathbf{p},c)= \left(\bm A_{i}^{(c)}\right)^{-1}(\mathbf{p},:)\bm b_{i}^{(c)}.
\end{equation}
Although the matrix $\bm A_{i}^{(c)}$ is sparse, its inverse is a dense matrix. As a result, $\delta \bm V_{i}(\mathbf{p},c)$ cannot be computed efficiently because computing the inverse of $\bm A_{i}^{(c)}$ for every column is effectively equivalent to solving the FOM. This issue arises from the fact that in a linear system of equations, the solution for any entry in the unknown vector is coupled to all other entries in that vector. The challenge is to solve for a select few entries without incurring the computational cost of solving the entire linear system.

To avoid this computational issue, we propose an algorithm that does not require inverting matrix $\bm A_{i}^{(c)}$. To this end, we note that the matrix $\delta \bm V_i = \bm V^k_{i+1} -\bm V^k_{i}$ can also be approximated accurately via a low-rank approximation since both  $\bm V^k_{i}$ and $\bm V^k_{i+1}$ are represented via low-rank approximations in the TDB-CUR.  In fact, if $\bm V^k_{i}$ and $\bm V^k_{i+1}$ are approximated with rank-$r$ matrices, then the maximum rank of $\delta \bm V_i$ is $2r$. 

We aim to construct a low-rank approximation for $\delta \bm{V}_i$. We begin by identifying a low-rank subspace for the columns of $\delta \bm{V}_i$. This subspace is spanned by the columns of $\mathbf{V}^{k-1}(:,\mathbf{s}^{k-1}) \in \mathbb{R}^{n\times r^{k-1}}$ and $\mathbf{V}^{k}(:,\mathbf{s}^{k}) \in \mathbb{R}^{n\times r^{k}}$
\begin{equation}\label{eq:V_delta}
\bm{V}_{\delta} = [\bm{V}^{k-1}(:,\bm{s}^{k-1}),\bm{V}^{k}(:,\bm{s}^{k})]\in \mathbb{R}^{n \times (r^{k-1}+r^k)},
\end{equation}
where $\mathbf{s}^{k-1}$ and $\mathbf{s}^{k}$ are the selected columns  and $r^{k-1}$ and $r^{k}$, are the ranks of the solution at time steps $k-1$ and $k$, respectively. 

The matrix $\bm{V}_{\delta}$ may contain singular values close to machine precision because $\bm{V}^{k-1}(:,\bm{s}^{k-1})$ and $\bm{V}^{k}(:,\bm{s}^{k})$ are highly correlated, being separated by only a small time step $\Delta t$. To construct a low-rank subspace for the columns of  $\bm{V}_{\delta}$ and remove the dimensions associated with small singular values we compute the SVD of $\bm{V}_{\delta}$, which is a thin matrix and its SVD can be computed efficiently:
\begin{equation}\label{eq:V_delta_SVD}
\bm{V}_{\delta}  \approx \bm{U}_{\delta}\bs{\Sigma}_{\delta}{\bm{Y}_{\delta}}^{\mathrm T},
\end{equation}
where  $\bm{U}_{\delta} \in \mathbb{R}^{n \times r_{\delta}}$ is the matrix of left singular values and $r_{\delta} \leq r^{k-1}+r^k$ is the rank of the approximation for $\delta \bm V_i$ and $\bs{\Sigma}_{\delta} = \mbox{diag}(\sigma_{\delta_1},\sigma_{\delta_2}, \dots, \sigma_{\delta_{r_\delta}})$. The matrix of right singular vectors $\bm{Y}_{\delta}$ is not used in the steps that follow. The truncation criterion retains only singular values $\sigma_{\delta_i}$ satisfying \cite{Penrose_1956}:
\begin{equation*}
\sigma_{\delta_i} > \epsilon_m \cdot \max(n,r^{k-1}+r^k) \cdot \sigma_{\delta_1}.
\end{equation*}

The orthonormal matrix $\bm{U}_{\delta}$ provide a basis for the columns of $\delta\bm{V}^k$.  By representing $\delta \bm V_{i}$ in a reduced basis $\bm{U}_{\delta}$, we can find $\delta \bm{V}_{i}(\mathbf{p},:)$ without forming or inverting $\bm A_{i}^{(c)}$. The correction matrix $\delta \bm{V}_{i}$, can be approximated with rank-$r_{\delta}$ approximation as follows:
\begin{equation}\label{eq:approx}
\delta \bm{V}_{i} \approx \widehat{\delta \bm{V}}_{i} =  \bm U_{\delta}\mathbf{Z}_{i},
\end{equation}
where $\mathbf{Z}_{i}\in \mathbb{R}^{r_{\delta}\times s}$ is the reduced coefficient matrix. 
Then, \cref{eq:approx} is substituted into \cref{eq:LSoE} to obtain the residual vector for each column as follows:
\begin{equation}\label{eq:res_delta}
    \bm R_{i}(:,c) = \bm{A}_{i}^{(c)}\bm{U}_{\delta}\mathbf{Z}_{i}(:,c) -  \bm{b}_{i}^{(c)},
\end{equation}
where $\bm R_{i}$ denotes the values of the residual at iteration $i$. 
The residual arises from the low-rank approximation error of $\delta \bm{V}_{i}$ in \cref{eq:approx}.
Therefore, the goal is to find $\mathbf{Z}_{i}(:,c)$ such that the residual $\bm R^{\delta}_{i}(:,c)$ is minimized in some sense. 

One approach to solving this problem is to find the least-squares solution, as shown below:
\begin{equation*}
\big(\bm{A}_{i}^{(c)}\bm{U}_{\delta}\big)^{\mathrm T}\big(\bm{A}_{i}^{(c)}\bm{U}_{\delta} \big)\mathbf{Z}_{i}(:,c) = \big(\bm{A}_{i}^{(c)}\bm{U}_{\delta}\big)^{\mathrm T}\bm{b}_{i}^{(c)}, \quad c=1,2, \dots, s.
\end{equation*}
The left-hand side now involves inverting a small matrix $\big(\bm{A}_{i}^{(c)}\bm{U}_{\delta}\big)^{\mathrm T}\big(\bm{A}_{i}^{(c)}\bm{U}_{\delta}\big) \in \mathbb{R}^{r_{\delta}\times r_{\delta}}$ instead of the large matrix $\bm{A}_i^{(c)}$, which can be solved with $\mathcal{O}(r^3)$ complexity, where we have used the fact that $r$ and $r_{\delta}$ are of the same order. However, this approach requires computing the action of matrix $\bm{A}_i^{(c)}$ on the columns of $\bm{U}$ and the vector $\bm{b}_i^{(c)}$ for every $c = 1, 2, \dots, s$. Even for sparse matrices $\bm{A}_i^{(c)}$, the cost of this operation scales with $\mathcal{O}(rsn)$, which is prohibitive for large $n$ and $s$.

To mitigate this issue, we use a residual collocation approach to find  $\mathbf{Z}_{i}(:,c)$ by setting the residual to zero at  DEIM-selected entries. To this end, we first compute the DEIM collocation points using the $\bm U_{\delta}$ as the basis:
\begin{equation}
    \bm p  = \texttt{DEIM}(\bm U_{\delta}),
\end{equation}
where $\bm p = [p_{1},p_{2}, \dots, p_{{r_{\delta}}}]$ is the row indices. Setting $\bm R_{i}(\bm p,c) = \bm 0$ results in:
\begin{equation}\label{eq:z}
\mathbf{Z}_{i}(:,c) = ( \underbrace{\bm{A}(\mathbf{p},:)_{i}^{(c)}\bm{U}_{\delta}}_{\bm A^{(c)}_r} )^{-1} \bm{b}_{i}^{(c)}(\mathbf{p}),
\end{equation}
where $\mathbf{A}^{(c)}_{r}$ is a reduced $r_{\delta}\times r_{\delta}$ matrix, and this smaller system can now be solved efficiently. It is worth noting that, by leveraging the sparsity inherent to the spatial discretization method, the large matrix $\bm{A}(\mathbf{p},:)_{i}^{(c)}$ can be multiplied by the vector $\bm{U}_{\delta}$ efficiently despite its size. The computation scales as $\mathcal{O}(sr^3)$, which is significantly smaller than the least-squares approach that scales as $\mathcal{O}(rsn)$. This reduction is achieved because the action of the matrix $\bm{A}_{i}^{(c)}$ on $\bm{U}_{\delta}$ does not need to be computed. Instead, only the action of $\mathcal{O}(r)$ rows of $\bm{A}_{i}^{(c)}$ on $\bm{U}_{\delta}$ is required. Another advantage of the above algorithm is that computing $\mathbf{Z}_{i}(:,c)$ according to \cref{eq:z} can be done independently for each $c$ index. As a result, the computation of $\mathbf{Z}_{i}(:,c)$ is highly parallelizable, similar to the calculation of the columns of $\bm{V}^k$.

Once the matrix $\bm Z_i$ is calculated, the low-rank construction of $\delta \bm V_i$ according to \cref{eq:approx} is complete. The Newton's correction matrix can now be evaluated efficiently at $\bm p$ rows:

\begin{equation}
    \widehat{\delta \bm V}_i(\bm p,:) = \bm U_{\delta}(\mathbf{p},:)\mathbf{Z}_{i}.
\end{equation}
Therefore, $\bm V^k(\bm p,:)$ can now be solved using Newton's method by iteratively  updating the solution at the selected rows:
\begin{equation}\label{eq:Newton_rows}
\bm V_{i+1}^{k}(\bm p,:) = \bm V_{i}^{k}(\bm p,:) + \bm U_{\delta}(\mathbf{p},:)\mathbf{Z}_{i},
\end{equation}
After updating the solution at selected rows, $\bm{A}(\mathbf{p},:)_{i}^{(c)}$ and $\bm{b}_{i}^{(c)}(\bm p)$  must be updated. However, updating these two quantities requires values of $\bm V_{i+1}^{k}$ at additional rows due to row dependency in the MDE. This is explained below in more detail.

The rows of MDEs obtained from the spatial discretization of parametric PDEs have dependencies. This means that advancing the state of any row to the next time step requires the values of other rows. The row dependency depends on the specific spatial discretization scheme, e.g., finite difference/element schemes or dense spectral methods. To advance the solution at DEIM-selected rows at each iteration of the nonlinear Newton solver, we need to update $\bm{A}(\mathbf{p},:)_{i}^{(c)}$ and $\bm{b}_{i}^{(c)}(\mathbf{p})$ based on the nonlinear state obtained from the previous iteration.  This includes the evaluation of $\bm F_i^k = \mathcal{F}(\hat{\bm V}_{i}^{k})$ and its Jacobian at rows with $\bm p$ indices.

Evaluating $\bm F_i^k(\bm p,:)$  requires the $\bm V(\bm p,:)$ as well as $\bm V(\bm p_a,:)$, where $\bm p_a$ is the indices of other rows whose values are required to evaluate $\bm F_i^k(\bm p,:)$. For example, let $f(v;x,t,\xi) = \partial^2 v/\partial x^2$ be the right-hand side of a one-dimensional PDE and $\bm F = \mathcal F(\bm V) = \bm D \bm V$, where $\bm D \in \mathbb{R}^{n \times n}$ is the discrete representation of $\partial^2 ( \cdot)/\partial x^2$ obtained via the second-order finite difference:
\begin{equation*}
    \bm F(i,:) = \frac{\bm V(i+1,:) -2\bm V(i,:)+\bm V(i-1,:)}{\Delta x^2}.
\end{equation*}
  Let $\bm p=[5]$ for simplicity. In this example, computing $\bm F(5,:)$ requires the values of $\{\bm V(4,:),\bm V(5,:),\bm V(6,:)\} $ and therefore, $\bm p_a=[4,6]$.  requires the state of the matrix at the dependent rows. Specifically, the spatial discretization requires the values of the set of adjacent rows, whose indices are denoted with $\bm p_a$.  However,  $\bm V_{i}^{k}(\mathbf{p_a},:)$ is not known. To resolve this issue, we leverage the correction basis $\bm U_{\delta}$ and coefficient $\mathbf{Z}_{i}(:,c)$ as a  low-rank estimation of our solution at the $\bm p_a$ rows:
  \begin{equation}\label{eq:Newton_rows_aux}
\bm V_{i+1}^{k}(\bm p_a,:) = \bm V_{i}^{k}(\bm p_a,:) + \bm U_{\delta}(\mathbf{p}_a,:)\mathbf{Z}_{i},
\end{equation}

Therefore, each  Newton's iteration for solving for selected rows involves updating $\bm{A}(\mathbf{p},:)_{i}^{(c)}$ and  $\bm{b}_{i}^{(c)}(\mathbf{p})$ using the values of $ \bm V_{i}^{k}(\bm p,:)$ and $ \bm V_{i}^{k}(\bm p_a,:)$. Next,  $\bm Z_i$ is computed via \cref{eq:z}. Then the solution at row induced of $\bm p$ and $\bm p_a$ is updated according to \cref{eq:Newton_rows} and \cref{eq:Newton_rows_aux}, respectively. These iterations continue until convergence is achieved. As iterations proceed, $\bm Z_i$ converges to zero, however, when oversampling is used the entries of $\bm Z_i$ converge to small non-zero values. This is explained in more detail in \cref{sec:os}. 

The implicit TDB-CUR algorithm is outlined in \cref{alg:ITDB}. The solution at $r$ selected columns is then computed implicitly using Newton's method. The rows are solved using the CUR Newton's method explained in this section. The stable CUR \cref{alg:SCUR} is then used to construct a rank-$r$ approximation using $\bm V^k(\bm p,:)$ and $\bm V^k(:,\bm s)$.  The solution matrix is then stored in the SVD-like factorized form: $\hat{\bm V}^k = \bm U^k \boldsymbol \Sigma^k \bm Y^{k^{\mathrm T}}$.  

\begin{algorithm}[htbp]
\caption{Implicit Time-Integration on Low-Rank Matrix Manifolds}
\hspace*{\algorithmicindent} \textbf{Input}: $\mathbf{U}^{k-1} \in \mathbb{R}^{n \times r}$, $\boldsymbol \Sigma^{k-1} \in \mathbb{R}^{r \times r}$, $\mathbf{Y}^{k-1} \in \mathbb{R}^{s \times r}$ \\
\hspace*{\algorithmicindent} \textbf{Output}: $\mathbf{U}^{k}$, $\boldsymbol \Sigma^{k}$, $\mathbf{Y}^{k}$
\begin{algorithmic}[1]
\State ${\bm s} \gets$ \texttt{DEIM}($\mathbf{Y}^{k-1}$)
\Comment{Compute $r$ columns indices}
\State $\mathbf{V}^{k-1}(:, \mathbf{s}) = \mathbf{U}^{k-1} \boldsymbol \Sigma^{k-1} {\mathbf{Y}(\mathbf{s},:)^{k-1}}^{\mathrm T} $
\Comment{Compute previous state at selected columns $\mathbf{s}$}
\State $\mathbf{V}^{k}(:, \mathbf{s}) = \texttt{Solve\_FOM}(\mathbf{V}^{k-1}(:, \mathbf{s}))$
\Comment{Advance} the selected columns
\State $\bm U_{\delta}, \bs{\Sigma}_{\delta}, \bm {Y}_{\delta} = \texttt{SVD}([\bm V^{k-1}(:, \bm s), \bm V^k(:, \bm s)])$ \Comment{Compute and truncate the orthonormal basis $\bm U_{\delta}$}
\State ${\bm p} \gets$ \texttt{DEIM-OS}($\bm U_{\delta},r_{\delta}+e$)
\Comment{Compute $r_{\delta}+e$ row indices}
\State $\bm p_a\gets$ \texttt{find\_adjacent}($\bm p$)\Comment{Find adjacent points required to compute $\bm V^k(\bm p,:)$}
\State $\mathbf{V}^{k-1}(\mathbf{p}, :) = \mathbf{U}^{k-1}(\mathbf{p}, :) \boldsymbol \Sigma^{k-1} {\mathbf{Y}^{k-1}}^{\mathrm T} $
\Comment{Compute previous state at selected rows $\mathbf{p}$}
\State $\mathbf{V}^{k-1}(\mathbf{p_a}, :) = \mathbf{U}^{k-1}(\mathbf{p_a}, :) \boldsymbol \Sigma^{k-1} {\mathbf{Y}^{k-1}}^{\mathrm T} $
\Comment{Compute previous state at selected rows $\mathbf{p_a}$}
\For {$c=1$ to $s$} \Comment{Loop over columns}
\State $i=1$ and $\epsilon = 1$
\State $\bm v_{1}^{k}(\bm p) = \bm V^{k-1}(\bm p,c)$  \Comment{Initialize the solution at $\bm p$}
\State $\bm v_{1}^{k}(\bm p_a) = \bm V^{k-1}(\bm p_a,c)$ \Comment{Initialize the solution at $\bm p_a$}
\While {$\epsilon > \epsilon_t$ \& $i \leq 6$} \Comment{Check the convergence and cut-off criteria}
\State \texttt{update}($\bm{A}_{i}(\bm p, :), \bm{b}_{i}(\bm p)$)
\Comment{Update $\bm{A}_{i}(\bm p, :), \bm{b}_{i}(\bm p)$ }
\State $\mathbf{z}^{(c)}= ( \bm{A}_i(\mathbf{p},:)^{(c)} \bm U_{\delta} )^{\dagger} \bm{b}_i^{(c)}(\mathbf{p})$ \Comment{Compute the reduced coordinates $\mathbf{z}^{(c)}$}
\State $\bm v_{i+1}^{k}(\bm p) = \bm v_{i}^{k}(\bm p) + \bm U_{\delta}(\mathbf{p},:) \mathbf{z}^{(c)}$ \Comment{Update the solution at $\bm p$}
\State $\bm v_{i+1}^{k}(\bm p_a) = \bm v_{i}^{k}(\bm p_a) + \bm U_{\delta}(\mathbf{p_a},:) \mathbf{z}^{(c)}$ \Comment{Update the solution at $\bm p_a$}
\State $\epsilon = \lVert \bm v_{i+1}^{k}(\bm p) - \bm v_{i}^{k}(\bm p) \rVert \times p^{-1}$ \Comment{Compute $\epsilon$ ($p=r_{\delta}+e$)}
\State $i=i+1$
\EndWhile
\State $\bm V^k(\bm p,c) = \bm v_{i+1}^{k}(\bm p)$ \Comment{Store the solution}
\EndFor
\State $\mathbf{U}^{k}, \boldsymbol \Sigma^{k}, \mathbf{Y}^{k} = \texttt{Stable\_CUR}(\mathbf{V}^{k}(:, \mathbf{s}), \bm V^k(\bm p,:))$
\Comment{Perform the Stable CUR algorithm}
\end{algorithmic}
\label{alg:ITDB}
\end{algorithm}

\subsubsection{CUR for Multistep Implicit Methods}\label{sec:multi-step}

Implicit multistep methods are widely employed for the time integration of stiff PDEs. These methods utilize information from the previous $l$ steps to calculate the solution at the current time step, with the number of prior steps determining the order of accuracy. Applying an implicit multistep time integration method to the time integration of the FOM results in:
\begin{equation*}
\sum_{j=0}^l a_j \bm V^{k-j}=\Delta t \sum_{j=0}^l b_j \mathcal{F}\left(\bm V^{k-j}\right),
\end{equation*}
with $a_l=1$. The coefficients $a_0, \ldots, a_{l-1}$ and $b_0, \ldots, b_l$ determine the method.  Replacing $\bm V^{k-j}$ with a rank-$r$ approximation solution denoted $\hat{\bm V}^{k-j}$  results in the following residual:
\begin{equation}\label{eq:Res_MS}
\bm R = \sum_{j=0}^l a_j \hat{\bm V}^{k-j}-\Delta t \sum_{j=0}^l b_j \mathcal{F}\left(\hat{\bm V}^{k-j}\right).
\end{equation}
Similar to the implicit Euler method, CUR Newton's method is employed to solve the nonlinear problem by setting \(\bm{R}(:,\bm{s}) = \bm{0}\) and \(\bm{R}(\bm{p},:)=\bm{0}\). The primary difference between the implicit Euler time integration and multistep schemes lies in the different expressions for the residuals, as given by \cref{eq:Residual} and \cref{eq:Res_MS}. Newton's method is applied to \cref{eq:Res_MS}, which involves solving a linear system of equations akin to \cref{eq:Newton5}, with the exception that \(\bm{A}_i^{(c)}\) and \(\bm{b}_i^{(c)}\) are derived for the multistep integration residual, as defined by \cref{eq:Res_MS}.

\subsubsection{CUR for Diagonally Implicit Runge-Kutta Methods}\label{sec:RK}

In this section, we present a CUR methodology for the time integration of the MDE given by \cref{eq:FOM1} using implicit Runge-Kutta (IRK) methods. We specifically focus on diagonally implicit Runge-Kutta (DIRK) methods, which are widely used due to their structure; in DIRK methods, each stage depends only on the current and previous stages, significantly simplifying both computation and implementation \cite{KC19}.   An $L$-stage DIRK method has the form:
\begin{equation*}
\bm V^{k} = \bm V^{k-1} + \Delta t \sum_{l=1}^L b_l \bm K^l ,
\end{equation*}
where,
\begin{equation}\label{eq:RK_st}
\bm K^{l}= \mathcal{F}\left(\bm V^{k-1}+\Delta t \sum_{m=1}^l a_{l m} \bm K^{m}\right), \quad l=1, \ldots, L.
\end{equation}
The coefficients $a_{lm}$ and $b_l$  define the Runge-Kutta and can be found in the Butcher tableau.

The above residual minimization problem is formulated using \(\bm{K}^l\) as the unknown, rather than \(\bm{V}\). It is also possible to formulate the DIRK versus the $l$-th stage values of \(\bm{V}\).   To achieve this, a low-rank approximation of \(\bm{K}^l\), denoted by \(\hat{\bm{K}}^l\), replaces \(\bm{K}^l\) in \cref{eq:RK_st}, which results in a residual as follows:
\begin{equation}\label{eq:RK_res}
\bm R = \hat{\bm K}^{l}- \mathcal{F}\left(\hat{\bm V}^{k-1}+\Delta t \sum_{m=1}^l a_{l m} \hat{\bm K}^{l}\right), \quad l=1, \ldots, L.
\end{equation}
The above residual minimization is defined as finding $ \hat{\bm K}^{l}$ such that $\bm R$ is minimized for each stage of the DIRK. Therefore, the DIRK consists of solving $L$ residual minimization problems.  
 We follow the CUR methodology used to solve for \(\hat{\bm V}^k\) in the multistep methods to solve for \(\hat{\bm{K}}^l\).

The linearized equations for solving \cref{eq:RK_res} can be derived as:
\begin{equation}\label{eq:LDIRK}
\bm A_{i}^{(c)} \delta \bm K^{l}_{i}(:,c) =\bm b^{(c)}_i, \quad c=1, 2, \dots, s, \quad l=1, \ldots, L,
\end{equation}
where
\begin{align}
\bm A_i^{(c)}  &= \bm I - a_{l l}\Delta t \mathcal{J}(\hat{\bm V}^{k-1}(:,c)+\Delta t \sum_{m=1}^{l-1} a_{l m} \hat{\bm K}^{m}(:,c) + \Delta t a_{l l} \hat{\bm K}_{i}^{l}(:,c)), \\
\bm b^{(c)}_i &=\mathcal{F}\left(\hat{\bm V}^{k-1}(:,c)+\Delta t \sum_{m=1}^{l-1} a_{l m} \hat{\bm K}^{m}(:,c) + \Delta t a_{l l} \hat{\bm K}_{i}^{l}(:,c)\right) - \hat{\bm K}_{i}^{l}(:,c).
\end{align}

We use one correction basis $\bm{U}_{\delta}$ for all stages of the DIRK. Therefore, the SVD should be performed using $\bm{V}^{k}(:,\bm{s})$ along with all stage solutions $\bm{K}^l(:, \bm{s})$:
\begin{equation*} 
\bm{V}_{\delta} = [\bm V^{k}(:, \bm s), \bm K^1(:, \bm s), \bm K^2(:, \bm s), \dots \bm K^L(:, \bm s)] \in \mathbb{R}^{n \times (L+1)r}.
\end{equation*}
Unlike single-stage methods, DIRK evaluates \(\mathcal{F}(\hat{\bm{V}})\) at multiple intermediate stages. Therefore, the correction space should include information from all stage solutions \(\bm{K}^l\), not just the final stage. This enriches the correction subspace, enabling accurate approximation of the solution at all DIRK stages using \(\bm{U}_{\delta}\) as the basis.

\subsection{Oversampling and Rank Adaptivity}
In this section, we introduce two modifications to the TDB-CUR algorithm: oversampling and rank adaptivity. Row oversampling enhances the accuracy of the CUR algorithm, while the rank adaptivity method adjusts the rank over time to control the low-rank approximation error.  
\subsubsection{Oversampling for Improved Condition Number}\label{sec:os}

Oversampling the rows, i.e. selecting more than $r_{\delta}$ rows, improves the condition number of DEIM algorithm by reducing $\eta_r = \|\bm U_{\delta}(\bm p,:)^{\dagger} \|$. We denote the number of oversampled rows with $e$. Therefore,  $\bm p$ contains $r_{\delta}+e$ row indices.  Oversampling requires a very minor modification to the TDB-CUR algorithm, in which the coefficient $\bm Z_i$ is computed as follows:
\begin{equation}
    \mathbf{Z}_{i}(:,c) = ( \underbrace{\bm{A}(\mathbf{p},:)_{i}^{(c)}\bm{U}_{\delta}}_{\bm A^{(c)}_r} )^{\dagger} \bm{b}_{i}^{(c)}(\mathbf{p}), \quad c=1,2,\dots, s.
\end{equation}
The number of DEIM-selected row indices is equal to \(r_{\delta}\). To determine the additional \(e\) row indices for oversampling, we employ the \texttt{GappyPOD+E} method \cite{peherstorfer2020stability}. Stated simply, oversampling changes the interpolation problem into a regression problem. Consequently, the residual at the selected rows can no longer be reduced to machine precision because the \(\bm{Z}\) coefficients are obtained by solving an overdetermined regression problem. As a result, the residual at the selected rows saturates to a small error, which is attributed to the low-rank approximation error.
   
\subsubsection{Rank Adaptivity}
To control the approximation error while minimizing computations, we employ the same adaptive rank criteria used in the explicit TDB-CUR method \cite{DPNFB23}. The rank is adjusted based on an error proxy \(\hat{\epsilon}\), which estimates the low-rank approximation error by computing the ratio singular value \(\hat{\sigma}_r\) to  the  Frobenius norm of $\hat{\bm V}$ as follows:
    \begin{equation*}
    \hat{\epsilon} = \frac{\hat{\sigma}_r(t)}{\left(\sum_{i=1}^r\hat{\sigma}_i(t)^2\right)^{1/2}}.
    \end{equation*}
    Rather than a fixed threshold, the rank is adjusted to keep $\hat{\epsilon}$ within a range $\epsilon_l \leq \hat{\epsilon} \leq \epsilon_u$ specified by the user. If $\hat{\epsilon} > \epsilon_u$, the rank is increased to $r+1$ to improve accuracy. If $\hat{\epsilon} < \epsilon_l$, the rank is decreased to $r-1$ to improve efficiency. This approach prevents excessive rank addition or removal that might occur with a hard threshold.

\section{\label{sec:DC}Demonstration cases}

\subsection{Advection-Diffusion Equation}

As the first example, we consider a one-dimensional stochastic advection-diffusion equation. Discretizing the diffusion term is a primary source of stiffness in PDEs with second-order derivative terms, as increasing spatial resolution introduces small time scales that impose stringent limits on \(\Delta t\) when using explicit schemes.   

We consider:
\begin{equation*}
\begin{aligned}
\frac{\partial v}{\partial t}+\frac{\partial v}{\partial x}=\alpha \frac{\partial^2 v}{\partial x^2}, && x \in[0,1], \ t \in[0,3],
\end{aligned}
\end{equation*}
where $v(x,t;\xi)$ is the random solution and  $\alpha = 0.1$. We consider Dirichlet boundary conditions at both boundaries, $v(0,t;\xi)=0$ and $v(1,t;\xi)=0$. We use second-order finite difference for the discretization of the first and second-order spatial derivatives. 
The grid consists of $n=201$ equidistant points the and random space is considered 20 dimensional in this problem ($d = 20$) and the number of samples is equal to $s=32$, which is sufficiently large for the purpose of this demonstration. To generate different samples, we start with a set of initial conditions that are parameterized by random variables. The initial condition $v(x,0;\xi)$ is given by a sum of Gaussian functions centered at $x_{0_i}$ with a fixed width:
\begin{equation*}
v(x,0;\xi) = \sum_{i=1}^{d} \exp\left(-\left(\frac{x - x_{0_i}}{0.5}\right)^2\right) \left(x(1-x)\right) \xi_i,
\end{equation*}
where $x_{0_i}$ are equidistant points in the interval $[0,1]$. The random variables $\xi_i$ are sampled from a standard normal distribution $\xi_i \sim \mathcal{N}(0,1)$.  The number of row oversampling is  $e=$15 in all the simulations. For the given parameters, the 4\textsuperscript{th} order explicit Runge-Kutta method is stable when the maximum time step size of $\Delta t=2\times 10^{-4}$. We solve the MDE using various implicit schemes with $\Delta t$ as large as $\Delta t= 0.5$.
Since the resulting MDE is linear, no Newton iteration is required and a single linear solve at each time step is performed to advance the solution implicitly. 

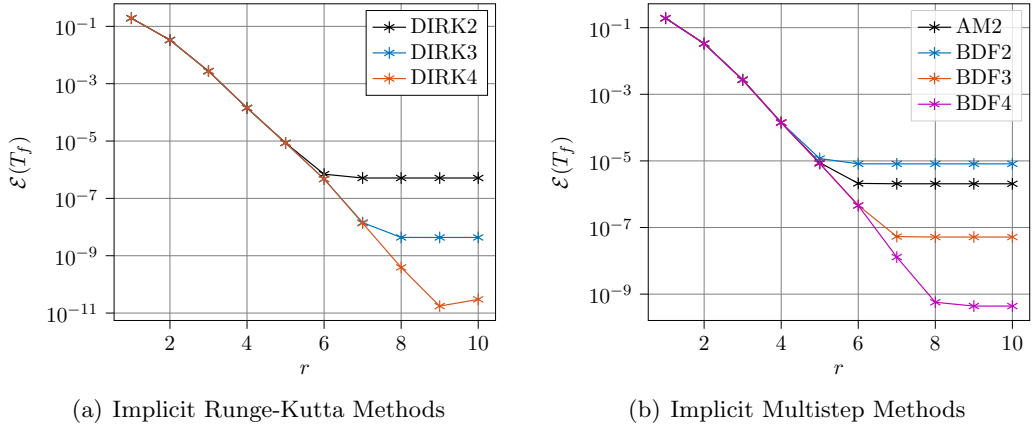
\begin{figure}
     \centering
     \subfigure[Implicit Runge-Kutta Methods]{
         \centering
\begin{tikzpicture}[scale=0.74]
\definecolor{gray}{RGB}{128,128,128}
\definecolor{myblue}{rgb}{0.00000,0.44700,0.74100}%
\definecolor{myred}{rgb}{0.85000,0.32500,0.09800}%
\definecolor{mygrey}{rgb}{.7 .7 .7}
\definecolor{myyellow}{rgb}{0.93 0.69 0.13}
\begin{axis}[
log basis y={10},
tick align=outside,
tick pos=left,
x grid style={gray},
xlabel={\(\displaystyle r\)},
xmajorgrids,
xmin=0.55, xmax=10.45,
xtick style={color=black},
y grid style={gray},
ylabel={\(\displaystyle \mathcal{E}(T_f)\)},
ymajorgrids,
ymin=5.51737863217086e-12, ymax=0.612193173487528,
ymode=log,
ytick style={color=black}
]
\addplot [semithick, black, mark=asterisk, mark size=3, mark options={solid}]
table {%
1 0.19267969534855
2 0.0334158237974237
3 0.00274072227188684
4 0.000142557515315242
5 8.61339611426419e-06
6 6.98496971211195e-07
7 5.14207479927522e-07
8 5.14028157723569e-07
9 5.14029763642634e-07
10 5.1403401941094e-07
};
\addlegendentry{DIRK2}
\addplot [semithick, myblue, mark=asterisk, mark size=3, mark options={solid}]
table {%
1 0.192679695347919
2 0.0334158237920603
3 0.00274072222354393
4 0.000142556586308064
5 8.59804825343358e-06
6 4.72958346937227e-07
7 1.42901353306587e-08
8 4.34731297645117e-09
9 4.33448876079052e-09
10 4.33863764407052e-09
};
\addlegendentry{DIRK3}
\addplot [semithick, myred, mark=asterisk, mark size=3, mark options={solid}]
table {%
1 0.192679695347917
2 0.0334158237920383
3 0.00274072222350904
4 0.000142556586268833
5 8.59804637055307e-06
6 4.72938426774958e-07
7 1.36195327753405e-08
8 3.90725298096431e-10
9 1.75301374026507e-11
10 2.99972686180317e-11
};
\addlegendentry{DIRK4}
\end{axis}

\end{tikzpicture}
         \label{Analytical_E_r-a}
     }
     \subfigure[Implicit Multistep Methods]{
         \centering
\begin{tikzpicture}[scale=0.74]
\definecolor{darkviolet1910191}{RGB}{191,0,191}
\definecolor{gray}{RGB}{128,128,128}
\definecolor{lightgray204}{RGB}{204,204,204}
\definecolor{myblue}{rgb}{0.00000,0.44700,0.74100}%
\definecolor{myred}{rgb}{0.85000,0.32500,0.09800}%
\definecolor{mygrey}{rgb}{.7 .7 .7}
\definecolor{myyellow}{rgb}{0.93 0.69 0.13}
\begin{axis}[
legend cell align={left},
legend style={fill opacity=0.8, draw opacity=1, text opacity=1, draw=lightgray204},
log basis y={10},
tick align=outside,
tick pos=left,
x grid style={gray},
xlabel={\(\displaystyle r\)},
xmajorgrids,
xmin=0.55, xmax=10.45,
xtick style={color=black},
y grid style={gray},
ylabel={\(\displaystyle \mathcal{E}(T_f)\)},
ymajorgrids,
ymin=1.61784659984662e-10, ymax=0.5212206327227,
ymode=log,
ytick style={color=black}
]
\addplot [semithick, black, mark=asterisk, mark size=3, mark options={solid}]
table {%
1 0.19267969535437
2 0.0334158238503385
3 0.00274072299324427
4 0.000142571413167741
5 8.84047531461951e-06
6 2.10980810210954e-06
7 2.05616282316921e-06
8 2.0561178122334e-06
9 2.05611779950153e-06
10 2.05611782298267e-06
};
\addlegendentry{AM2}
\addplot [semithick, myblue, mark=asterisk, mark size=3, mark options={solid}]
table {%
1 0.192220229196348
2 0.0333541000429775
3 0.00273223202618956
4 0.000141865443092943
5 1.17884273469511e-05
6 8.16701268314515e-06
7 8.15376155788236e-06
8 8.15375045289252e-06
9 8.15375044573102e-06
10 8.15375049640632e-06
};
\addlegendentry{BDF2}
\addplot [semithick, myred, mark=asterisk, mark size=3, mark options={solid}]
table {%
1 0.192217182522468
2 0.033294232455113
3 0.00272124870668837
4 0.0001407992886569
5 8.44415799707085e-06
6 4.63081211530363e-07
7 5.33999107550163e-08
8 5.17629030351701e-08
9 5.1761966570776e-08
10 5.17619182439573e-08
};
\addlegendentry{BDF3}
\addplot [semithick, darkviolet1910191, mark=asterisk, mark size=3, mark options={solid}]
table {%
1 0.192215288214167
2 0.0332341341259263
3 0.00271032133220648
4 0.00013997222166029
5 8.37495134596736e-06
6 4.55270879156044e-07
7 1.29485677760357e-08
8 5.70226486331095e-10
9 4.37712229177192e-10
10 4.37646025373581e-10
};
\addlegendentry{BDF4}
\end{axis}

\end{tikzpicture}
         \label{Analytical_E_r-b}
         }
        \caption{Error at the final time $(T_f=3)$ versus rank for different methods.}
        \label{Analytical_E_r}
\end{figure}

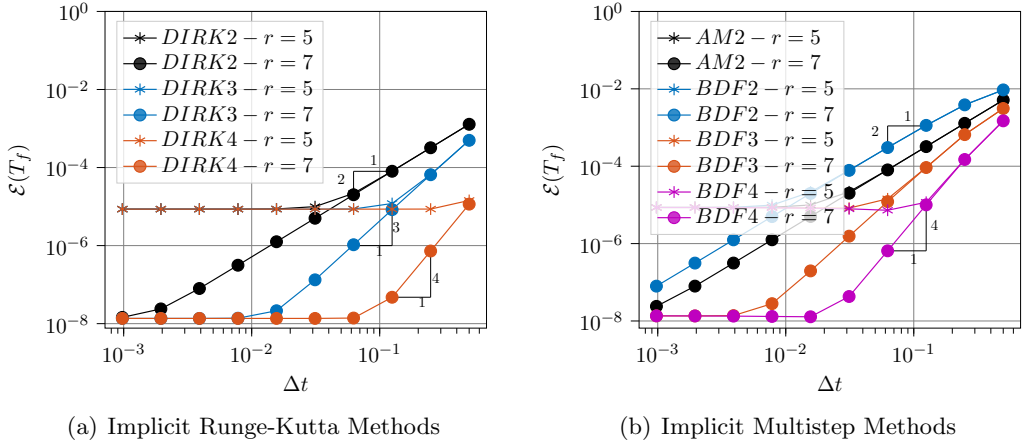
\begin{figure}
     \centering
     \subfigure[Implicit Runge-Kutta Methods]{
         \centering
\begin{tikzpicture}[scale=0.74]
\definecolor{gray}{RGB}{128,128,128}
\definecolor{lightgray204}{RGB}{204,204,204}
\definecolor{myblue}{rgb}{0.00000,0.44700,0.74100}%
\definecolor{myred}{rgb}{0.85000,0.32500,0.09800}%
\definecolor{mygrey}{rgb}{.7 .7 .7}
\definecolor{myyellow}{rgb}{0.93 0.69 0.13}
\begin{axis}[
legend cell align={left},
legend style={
  fill opacity=0.8,
  draw opacity=1,
  text opacity=1,
  at={(0.03,0.97)},
  anchor=north west,
  draw=lightgray204
},
log basis x={10},
log basis y={10},
tick align=outside,
tick pos=left,
x grid style={gray},
xlabel={\(\displaystyle \Delta t\)},
xmajorgrids,
xmin=0.000714885593723449, xmax=0.683020128377198,
xmode=log,
xtick style={color=black},
y grid style={gray},
ylabel={\(\displaystyle \mathcal{E}(T_f)\)},
ymajorgrids,
ymin=7.68147328720559e-09, ymax=1,
ymode=log,
ytick style={color=black}
]
\addplot [semithick, black, mark=asterisk, mark size=3, mark options={solid}]
table {%
0.5 0.00128367666493902
0.25 0.000321246720086493
0.125 8.07669222998925e-05
0.0625 2.18384033147742e-05
0.03125 9.95473385521272e-06
0.015625 8.68893174065155e-06
0.0078125 8.60376446910221e-06
0.00390625 8.5984042691924e-06
0.001953125 8.59806778341676e-06
0.0009765625 8.59804676052301e-06
};
\addlegendentry{$DIRK2 - r=5$}
\addplot [semithick, black, mark=*, mark size=3, mark options={solid}]
table {%
0.5 0.00128383710272449
0.25 0.000321193416540209
0.125 8.03125241643698e-05
0.0625 2.00790137942335e-05
0.03125 5.01982750512091e-06
0.015625 1.2550280040159e-06
0.0078125 3.14028363645771e-07
0.00390625 7.96070106338802e-08
0.001953125 2.38746890290871e-08
0.0009765625 1.44745262268683e-08
};
\addlegendentry{$DIRK2 - r=7$}
\addplot [semithick, myblue, mark=asterisk, mark size=3, mark options={solid}]
table {%
0.5 0.000498199451557648
0.25 6.6007654966537e-05
0.125 1.19764083980945e-05
0.0625 8.66179646630981e-06
0.03125 8.5990187825165e-06
0.015625 8.5980568690404e-06
0.0078125 8.59804581582282e-06
0.00390625 8.59804537417558e-06
0.001953125 8.59804536016198e-06
0.0009765625 8.59804535909605e-06
};
\addlegendentry{$DIRK3 - r=5$}
\addplot [semithick, myblue, mark=*, mark size=3, mark options={solid}]
table {%
0.5 0.000498264993676653
0.25 6.54730146617554e-05
0.125 8.34550377387937e-06
0.0625 1.05218444727907e-06
0.03125 1.32735210257767e-07
0.015625 2.14170407584134e-08
0.0078125 1.37744104993172e-08
0.00390625 1.36219193166548e-08
0.001953125 1.36196154001872e-08
0.0009765625 1.36195259694766e-08
};
\addlegendentry{$DIRK3 - r=7$}
\addplot [semithick, myred, mark=asterisk, mark size=3, mark options={solid}]
table {%
0.5 1.43566003463264e-05
0.25 8.6287961964916e-06
0.125 8.5982047883048e-06
0.0625 8.59805312644414e-06
0.03125 8.59804706622615e-06
0.015625 8.59804831954656e-06
0.0078125 8.59804548357067e-06
0.00390625 8.59804536714551e-06
0.001953125 8.59804535936447e-06
0.0009765625 8.59804536113438e-06
};
\addlegendentry{$DIRK4 - r=5$}
\addplot [semithick, myred, mark=*, mark size=3, mark options={solid}]
table {%
0.5 1.15432265805238e-05
0.25 7.26920915146299e-07
0.125 4.74537556629603e-08
0.0625 1.39124930198523e-08
0.03125 1.36207221472729e-08
0.015625 1.36195339053031e-08
0.0078125 1.36195286429379e-08
0.00390625 1.3619523849014e-08
0.001953125 1.36195448542219e-08
0.0009765625 1.36195259305577e-08
};
\addlegendentry{$DIRK4 - r=7$}
 \draw[black, thin] (axis cs:0.0625,0.00008) -- (axis cs:0.125,0.00008); 
 \draw[black, thin] (axis cs:0.0625,0.00008) -- (axis cs:0.0625,0.00002); 
 \node[scale=0.7] at (axis cs:0.05,0.00005) [anchor=center] {$2$};
 \node[scale=0.7] at (axis cs:0.09,0.00015) [anchor=center] {$1$};
 \draw[black, thin] (axis cs:0.125,0.000001) -- (axis cs:0.125,0.0000085); 
\draw[black, thin] (axis cs:0.0625,0.000001) -- (axis cs:0.125,0.000001); 
 \node[scale=0.7] at (axis cs:0.1,0.0000007) [anchor=center] {$1$};
 \node[scale=0.7] at (axis cs:0.135,0.000003) [anchor=center] {$3$};
 \draw[black, thin] (axis cs:0.25,0.0000000475) -- (axis cs:0.25,0.00000073); 
 \draw[black, thin] (axis cs:0.125,0.0000000475) -- (axis cs:0.25,0.0000000475); 
\node[scale=0.7] at (axis cs:0.215,0.000000035) [anchor=center] {$1$};
\node[scale=0.7] at (axis cs:0.275,0.00000015) [anchor=center] {$4$};
\end{axis}
\end{tikzpicture}
         \label{Analytical_E_dt-a}
     }
     \subfigure[Implicit Multistep Methods]{
         \centering
          \begin{tikzpicture}[scale=0.74]
\definecolor{darkviolet1910191}{RGB}{191,0,191}
\definecolor{gray}{RGB}{128,128,128}
\definecolor{lightgray204}{RGB}{204,204,204}
\definecolor{myblue}{rgb}{0.00000,0.44700,0.74100}
\definecolor{myred}{rgb}{0.85000,0.32500,0.09800}
\definecolor{mygrey}{rgb}{.7 .7 .7}
\definecolor{myyellow}{rgb}{0.93 0.69 0.13}

\begin{axis}[
    legend cell align={left},
    legend style={
      fill opacity=0.8,
      draw opacity=1,
      text opacity=1,
      at={(0.03,0.97)},
      anchor=north west,
      draw=lightgray204,
      scale = 0.5
    },
    log basis x={10},
    log basis y={10},
    tick align=outside,
    tick pos=left,
    x grid style={gray},
    xlabel={\(\displaystyle \Delta t\)},
    xmajorgrids,
    xmin=0.000714885593723449, xmax=0.683020128377198,
    xmode=log,
    xtick style={color=black},
    y grid style={gray},
    ylabel={\(\displaystyle \mathcal{E}(T_f)\)},
    ymajorgrids,
    ymin=6.53547073157355e-09, ymax=1,
    ymode=log,
    ytick style={color=black}
]
\addplot [semithick, black, mark=asterisk, mark size=3, mark options={solid}]
table {%
0.5 0.00511843972487392
0.25 0.00128386584130463
0.125 0.00032130846399206
0.0625 8.07714499440461e-05
0.03125 2.18424564621489e-05
0.015625 9.95614531807383e-06
0.0078125 8.68914797840515e-06
0.00390625 8.60376750216908e-06
0.001953125 8.5984030968303e-06
0.0009765625 8.59806771796775e-06
};
\addlegendentry{$AM2 - r=5$}
\addplot [semithick, black, mark=*, mark size=3, mark options={solid}]
table {%
0.5 0.00511843275926642
0.25 0.0012838371037823
0.125 0.000321193416745073
0.0625 8.03125234868716e-05
0.03125 2.00790123719864e-05
0.015625 5.01982401745246e-06
0.0078125 1.25502760993832e-06
0.00390625 3.14033103054471e-07
0.001953125 7.96074204311411e-08
0.0009765625 2.38741926477412e-08
};
\addlegendentry{$AM2 - r=7$}
\addplot [semithick, myblue, mark=asterisk, mark size=3, mark options={solid}]
table {%
0.5 0.00938764346879118
0.25 0.00385645210718292
0.125 0.00113667266216489
0.0625 0.000303559950630805
0.03125 7.85758227556622e-05
0.015625 2.15423769714025e-05
0.0078125 9.88212445134755e-06
0.00390625 8.65573984769826e-06
0.001953125 8.58714926610737e-06
0.0009765625 8.59008816895472e-06
};
\addlegendentry{$BDF2 - r=5$}
\addplot [semithick, myblue, mark=*, mark size=3, mark options={solid}]
table {%
0.5 0.00938783966014351
0.25 0.00385649545718367
0.125 0.00113665624150242
0.0625 0.000303453759242504
0.03125 7.81322530539849e-05
0.015625 1.98086336457052e-05
0.0078125 4.98615497176331e-06
0.00390625 1.25082645678409e-06
0.001953125 3.1350605695847e-07
0.0009765625 7.95386328921709e-08
};
\addlegendentry{$BDF2 - r=7$}
\addplot [semithick, myred, mark=asterisk, mark size=3, mark options={solid}]
table {%
0.5 0.00310642924339622
0.25 0.00065572139881684
0.125 9.28094267755856e-05
0.0625 1.43739870933106e-05
0.03125 8.26911603035213e-06
0.015625 8.36028642926129e-06
0.0078125 8.47763329563063e-06
0.00390625 8.53770992054791e-06
0.001953125 8.56785016150473e-06
0.0009765625 8.58294097578056e-06
};
\addlegendentry{$BDF3 - r=5$}
\addplot [semithick, myred, mark=*, mark size=3, mark options={solid}]
table {%
0.5 0.00310656221925954
0.25 0.000655713678734701
0.125 9.25594203861688e-05
0.0625 1.21641692548088e-05
0.03125 1.55617693989506e-06
0.015625 1.97104085042984e-07
0.0078125 2.80377358281035e-08
0.00390625 1.37778369017333e-08
0.001953125 1.35276037375513e-08
0.0009765625 1.35708084440119e-08
};
\addlegendentry{$BDF3 - r=7$}
\addplot [semithick, darkviolet1910191, mark=asterisk, mark size=3, mark options={solid}]
table {%
0.5 0.00148847592723507
0.25 0.00014924893545891
0.125 1.17086496208138e-05
0.0625 7.27583847196493e-06
0.03125 7.9091244679315e-06
0.015625 8.2505020202248e-06
0.0078125 8.42354300529318e-06
0.00390625 8.5106226933337e-06
0.001953125 8.55429151991364e-06
0.0009765625 8.57615787997132e-06
};
\addlegendentry{$BDF4 - r=5$}
\addplot [semithick, darkviolet1910191, mark=*, mark size=3, mark options={solid}]
table {%
0.5 0.00148881095739272
0.25 0.000149209320788956
0.125 1.00363166379234e-05
0.0625 6.50799284247256e-07
0.03125 4.29155336773623e-08
0.015625 1.28375190810462e-08
0.0078125 1.30883620383791e-08
0.00390625 1.33514479543319e-08
0.001953125 1.34849849548021e-08
0.0009765625 1.35521341199819e-08
};
\addlegendentry{$BDF4 - r=7$}

\draw[black, thin] (axis cs:0.0625,0.0011) -- (axis cs:0.125,0.0011); 
\draw[black, thin] (axis cs:0.0625,0.0003) -- (axis cs:0.0625,0.0011); 
\node[scale=0.7] at (axis cs:0.05,0.0008) [anchor=center] {$2$};
\node[scale=0.7] at (axis cs:0.09,0.002) [anchor=center] {$1$};


\draw[black, thin] (axis cs:0.125,0.00000065) -- (axis cs:0.125,0.00001); 
\draw[black, thin] (axis cs:0.0625,0.00000065) -- (axis cs:0.125,0.00000065); 
\node[scale=0.7] at (axis cs:0.145,0.000003) [anchor=center] {$4$};
\node[scale=0.7] at (axis cs:0.1,0.0000004) [anchor=center] {$1$};
\end{axis}
\end{tikzpicture}
    \label{Analytical_E_dt-b}
    }
    
        \caption{Error at the final time $(T_f=3)$ versus step size $\Delta t$  for different methods for r = 5, 7. }
        \label{Analytical_E_dt}
\end{figure}

We report on the convergence study of the proposed methodology for the following time integration schemes: second to fourth-order diagonally implicit Runge-Kutta (DIRK2-DIRK4) and multistep methods including the second-order Adams-Moulton (AM2) and the second to fourth-order backward differentiation formula (BDF2-BDF4).  To evolve the analytical solution of this linear system in time:
\begin{equation}\label{eq:adv_diff}
\frac{d\mathbf{V}}{d t}=\bm L \mathbf{V},
\end{equation}
the matrix exponential is used as
$\mathbf{V}(t)=e^{\bm L t} \mathbf{V}_0$, where $\bm V_0 \in \mathbb{R}^{n\times s}$ is the matrix of the random initial condition and $\bm L \in \mathbb{R}^{n\times n}$ is the finite-difference discrete representation of the first and second spatial derivatives.  The first and last rows of $\bm L$ are set to zero to enforce the homogeneous Dirichlet boundary conditions.  It is easy to show that the MDE given by \cref{eq:adv_diff} has an exact rank of $r=d=20$.

The error in each iteration equals the difference between the exact solution and the approximated solution based on the TDB.  The relative error is computed as follows:
\begin{equation}\label{eq:r_error}
\mathcal{E}(t) = \frac{\big \| \mathbf{V}_R^{k} - \mathbf{V}^{k}\big \|_F }{ \big \| \mathbf{V}_R^{k}\big \|_F}.
\end{equation}

The approximation of TDB-CUR to the MDE given by \cref{eq:adv_diff} contains two types of errors: the temporal error and the low-rank approximation error. In the following, we perform convergence studies for both of these errors. We first report the convergence with respect to varying rank. 
In \cref{Analytical_E_r-a}, the error at the final time ($T_f=3$) versus rank is shown for DIRK2, DIRK3, and DIRK4. All of these cases have the same time advancement of $\Delta t  = 0.01$. For smaller ranks, the low-rank approximation error is dominant, and therefore all three methods show the same error. As the rank increases, the temporal error becomes dominant, and the error of the high-order IRK methods saturates at smaller values.
In \cref{Analytical_E_r-b}, we observe a similar behavior for four different multistep methods.  

Next, we report on the temporal convergence. In \cref{Analytical_E_dt-a} we conduct a comparative convergence study for DIRK2, DIRK3, and DIRK4. We present the error as a function of the time-step size $\Delta t$ for two different ranks, $r = 5$ and $r = 7$. These results confirm that DIRK2, DIRK3, and DIRK4 achieve second-order, third-order, and fourth-order accuracy, respectively. At smaller values of $\Delta t$, the low-rank error becomes dominant. It is evident that the error saturation levels for $r = 7$ are lower than those for $r = 5$.  Also, the higher order DIRK methods achieve the saturation level at larger $\Delta t$s than the lower order ones.   A similar study is conducted for implicit multistep methods and the results are shown in \cref{Analytical_E_dt-b}. The multistep methods behave similarly to the DIRK methods.  

\subsection{Stochastic Burgers’ Equation}

The second example is a one-dimensional Burgers equation   
\begin{equation*}
\frac{\partial v}{\partial t}+v \frac{\partial v}{\partial x} = \nu \frac{\partial^{2} v}{\partial x^{2}}, \quad  x \in[0,1], t \in[0,1], 
\end{equation*}
subject to random initial conditions as follows:
\begin{equation*}
v(x, 0 ; \xi) = 0.5\sin(2\pi x)\left(e^{\cos(2\pi x)}-1.5\right)+\sigma_{x} \sum_{i=1}^{d} \sqrt{\lambda_{x_{i}}} \psi_{i}(x) \xi_{i}, \quad    \xi_{i}\sim \mathcal{N}(0,1),
\end{equation*}
where $\nu=0.01$ and Dirichlet boundary conditions are imposed at both boundaries where $v(0,t;\xi)=0$ and $v(1,t;\xi)=0$ . We consider a four-dimensional random space ($d = 4$) in this example and $\xi_{i}$'s are sampled from a normal distribution with $s= 32$. In the above equation, $\lambda_{x_{i}}$ and $\psi_{i}(x)$ are the eigenvalues and eigenvectors of the spatial kernel, $k$, with $\sigma_{x}=0.001$ \footnote{The code for the stochastic Burgers' equation example, along with the $\lambda_{x}$ and $\psi(x)$ data, can be accessed from: github.com/BabaeeLab/Implicit-TDB-CUR.}. The equation is solved numerically using the central second-order finite difference method on a uniform grid with $n=512$ grid points. To allow the dynamics to develop before deploying the low-rank approximation, the problem is first evolved up to $t=0.1$ using the full-order model. For all experiments, the DIRK4 method and the following parameters are used unless stated otherwise: $\epsilon_l=10^{-9}$, $\epsilon_u=10^{-8}$, $\epsilon_t=10^{-14}$, $e=15$, and $\Delta t =0.01$. For the mentioned parameters, the largest $\Delta t $ for which the explicit Runge–Kutta method is stable is roughly equal to $\Delta t =2 \times 10^{-4}$.

\begin{figure}
     \centering
     \subfigure[]{
         \centering
         \scalebox{.85}{\input{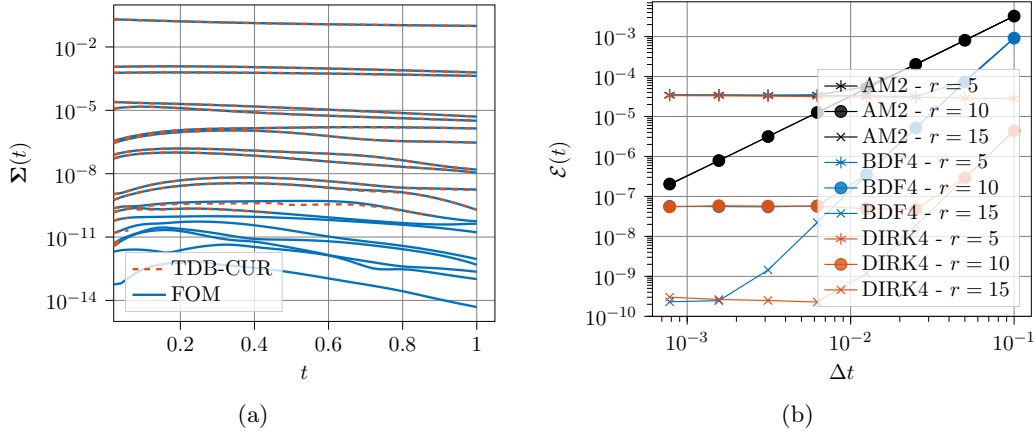}}
         \label{fig:burger_SC-a}
     }
     \subfigure[]{
         \centering
         \scalebox{1.05}{
\begin{tikzpicture}[scale=.7]

\definecolor{myblue}{rgb}{0.00000,0.44700,0.74100}%
\definecolor{myred}{rgb}{0.85000,0.32500,0.09800}%
\definecolor{myblack}{rgb}{0,0,0}%

\definecolor{gray}{RGB}{128,128,128}
\definecolor{lightgray204}{RGB}{204,204,204}

\begin{axis}[
legend cell align={left},
legend style={
  fill opacity=0.8,
  draw opacity=1,
  text opacity=1,
  at={(0.97,0.03)},
  anchor=south east,
  draw=lightgray204
},
log basis x={10},
log basis y={10},
tick align=outside,
tick pos=left,
x grid style={gray},
xlabel={\(\displaystyle \Delta t\)},
xmajorgrids,
xmin=0.000612956326481837, xmax=0.127456062731926,
xmode=log,
xtick style={color=black},
y grid style={gray},
ylabel={\(\displaystyle \mathcal{E}(t)\)},
ymajorgrids,
ymin=9.98929037754232e-11, ymax=0.00742151830933943,
ymode=log,
ytick style={color=black}
]
\addplot [semithick, myblack, mark=asterisk, mark size=3, mark options={solid}]
table {%
0.1 0.00325418967182974
0.05 0.000807068490026933
0.025 0.000199552269301474
0.0125 5.65011700184739e-05
0.00625 3.44421842189192e-05
0.003125 3.38604372098527e-05
0.0015625 3.43195858186479e-05
0.00078125 3.44115799584324e-05
};
\addlegendentry{AM2 - $r=5$}
\addplot [semithick, myblack, mark=*, mark size=3, mark options={solid}]
table {%
0.1 0.00325628142205457
0.05 0.000808467486355307
0.025 0.000201771898968316
0.0125 5.04217303034039e-05
0.00625 1.26034390510632e-05
0.003125 3.15083691948586e-06
0.0015625 7.89099922146943e-07
0.00078125 2.04226196219612e-07
};
\addlegendentry{AM2 - $r=10$}
\addplot [semithick, myblack, mark=x, mark size=3, mark options={solid}]
table {%
0.1 0.00325628148106808
0.05 0.000808468647470377
0.025 0.000201772945222889
0.0125 5.04218311659957e-05
0.00625 1.26041400836632e-05
0.003125 3.1509713364704e-06
0.0015625 7.87768754181057e-07
0.00078125 1.96975119887603e-07
};
\addlegendentry{AM2 - $r=15$}
\addplot [semithick, myblue, mark=asterisk, mark size=3, mark options={solid}]
table {%
0.1 0.000915931782721558
0.05 7.87283233600428e-05
0.025 3.12227196499065e-05
0.0125 3.30949317655405e-05
0.00625 3.37340810917031e-05
0.003125 3.39587690781572e-05
0.0015625 3.43653746668775e-05
0.00078125 3.44665101977189e-05
};
\addlegendentry{BDF4 - $r=5$}
\addplot [semithick, myblue, mark=*, mark size=3, mark options={solid}]
table {%
0.1 0.00091295949304541
0.05 7.33387969393379e-05
0.025 5.19160432432896e-06
0.0125 3.49956481157517e-07
0.00625 5.78129065773922e-08
0.003125 5.48940171949265e-08
0.0015625 5.5673893540117e-08
0.00078125 5.64309101471671e-08
};
\addlegendentry{BDF4 - $r=10$}
\addplot [semithick, myblue, mark=x, mark size=3, mark options={solid}]
table {%
0.1 0.000912960454684293
0.05 7.33375719343772e-05
0.025 5.19069702656126e-06
0.0125 3.4496590720029e-07
0.00625 2.2218418516658e-08
0.003125 1.44262597002707e-09
0.0015625 2.46575375810195e-10
0.00078125 2.31473466294555e-10
};
\addlegendentry{BDF4 - $r=15$}
\addplot [semithick, myred, mark=asterisk, mark size=3, mark options={solid}]
table {%
0.1 2.80859299726189e-05
0.05 2.88896494428629e-05
0.025 3.04767871189589e-05
0.0125 3.11625718960397e-05
0.00625 3.14135009857913e-05
0.003125 3.22695092939798e-05
0.0015625 3.29740135790979e-05
0.00078125 3.3328052843093e-05
};
\addlegendentry{DIRK4 - $r=5$}
\addplot [semithick, myred, mark=*, mark size=3, mark options={solid}]
table {%
0.1 4.38540642393421e-06
0.05 2.89713756788259e-07
0.025 4.77363100288788e-08
0.0125 4.93448929339989e-08
0.00625 5.69502199948653e-08
0.003125 5.72710789382183e-08
0.0015625 5.84219026558391e-08
0.00078125 5.49924926326615e-08
};
\addlegendentry{DIRK4 - $r=10$}
\addplot [semithick, myred, mark=x, mark size=3, mark options={solid}]
table {%
0.1 4.3852281986664e-06
0.05 2.87590483382466e-07
0.025 1.85435604996407e-08
0.0125 1.19450091590546e-09
0.00625 2.27669818672806e-10
0.003125 2.48837374242485e-10
0.0015625 2.63654492697006e-10
0.00078125 3.00771471212367e-10
};
\addlegendentry{DIRK4 - $r=15$}
\end{axis}

\end{tikzpicture}}
         \label{fig:burger_SC-b}
    }
        \caption{Burgers' Equation: Left panel shows the singular values over time for FOM and TDB-CUR. Right panel shows the relative error vs time-step size $\Delta t$ for AM2, BDF4, and IRK4 integrators at ranks $r=5,10,15$.}
        \label{fig:burger_SC}
\end{figure}

\cref{fig:burger_SC-a} shows the evolution of singular values over time for both the full order model (FOM) and the proposed implicit TDB-CUR method with DIRK4 as the time integration method. The FOM singular values are obtained by computing the SVD at each time step of the FOM solution.  It can be observed that TDB-CUR accurately matches the leading full-order singular values. Due to rank adaptivity, the number of singular values of TDB-CUR, which is equal to the rank of the approximation, varies with time. 
\cref{fig:burger_SC-b} examines the convergence of the method for different implicit time integrators - AM2, BDF4, and DIRK4.  The error versus time-step size $\Delta t$ is shown for reduced ranks $r=5, 10, 15$. As expected, BDF4 and DIRK4 display fourth-order convergence, while AM2 has second-order accuracy. For computing the relative error we use \cref{eq:r_error} and reference solution is calculated with the fourth-order explicit Runge-Kutta method using a time step size of $\Delta t = 10^{-7}$. For a fixed $\Delta t$, increasing the rank $r$ decreases the low-rank approximation error.  All integrators reach an error plateau at small $\Delta t$, corresponding to the optimal low-rank approximation error for the rank $r$ used. However, for larger $\Delta t$, the temporal integration error dominates over the total error.

\begin{figure}
     \centering
     \subfigure[]{
         \centering
         \scalebox{.75}{
\begin{tikzpicture}[scale=1]

\definecolor{myblue}{rgb}{0.00000,0.44700,0.74100}%
\definecolor{myred}{rgb}{0.85000,0.32500,0.09800}%
\definecolor{myblack}{rgb}{0,0,0}%

\definecolor{darkgray176}{RGB}{176,176,176}
\definecolor{lightgray204}{RGB}{204,204,204}

\begin{axis}[
legend cell align={left},
legend style={fill opacity=0.8, draw opacity=1, text opacity=1, draw=lightgray204},
tick align=outside,
tick pos=left,
x grid style={darkgray176},
xlabel={\(\displaystyle t\)},
xmin=-0.029, xmax=1.049,
xtick style={color=black},
y grid style={darkgray176},
ylabel={\(\displaystyle r(t)\)},
ymin=1, ymax=60,
ytick style={color=black}
]
\addplot [line width=1.2, black]
table {%
0.02 19
0.04 17
0.06 15
0.08 13
0.1 11
0.12 9
0.14 7
0.16 7
0.18 7
0.2 7
0.22 7
0.24 7
0.26 7
0.28 7
0.3 7
0.32 7
0.34 7
0.36 7
0.38 7
0.4 7
0.42 7
0.44 7
0.46 7
0.48 7
0.5 7
0.52 7
0.54 7
0.56 7
0.58 7
0.6 7
0.62 7
0.64 7
0.66 7
0.68 7
0.7 7
0.72 7
0.74 7
0.76 7
0.78 7
0.8 7
0.82 7
0.84 7
0.86 7
0.88 7
0.9 7
0.92 7
0.94 7
0.96 7
0.98 7
1 7
};
\addlegendentry{$r$ - $\epsilon_l=1 \times 10^{-6}$}
\addplot [line width=1.2, myblue]
table {%
0.02 33
0.04 26
0.06 27
0.08 26
0.1 22
0.12 18
0.14 14
0.16 14
0.18 14
0.2 14
0.22 14
0.24 14
0.26 14
0.28 14
0.3 14
0.32 14
0.34 14
0.36 14
0.38 14
0.4 14
0.42 14
0.44 14
0.46 14
0.48 14
0.5 14
0.52 14
0.54 14
0.56 14
0.58 14
0.6 14
0.62 14
0.64 14
0.66 14
0.68 14
0.7 14
0.72 14
0.74 14
0.76 14
0.78 14
0.8 14
0.82 14
0.84 14
0.86 14
0.88 14
0.9 14
0.92 14
0.94 14
0.96 14
0.98 14
1 14
};
\addlegendentry{$r_{\delta}$ - $\epsilon_l=1 \times 10^{-6}$}
\addplot [line width=1.2, myblack,mark=x]
table {%
0.02 19
0.04 17
0.06 15
0.08 13
0.1 11
0.12 9
0.14 9
0.16 9
0.18 9
0.2 9
0.22 9
0.24 9
0.26 9
0.28 9
0.3 9
0.32 9
0.34 9
0.36 9
0.38 9
0.4 9
0.42 9
0.44 9
0.46 9
0.48 9
0.5 9
0.52 9
0.54 9
0.56 9
0.58 9
0.6 9
0.62 9
0.64 9
0.66 9
0.68 9
0.7 9
0.72 9
0.74 9
0.76 9
0.78 9
0.8 9
0.82 9
0.84 9
0.86 9
0.88 9
0.9 9
0.92 9
0.94 9
0.96 9
0.98 9
1 9
};
\addlegendentry{$r$ - $\epsilon_l=1 \times 10^{-8}$}
\addplot [line width=1.2, myblue, mark=x]
table {%
0.02 33
0.04 26
0.06 27
0.08 26
0.1 22
0.12 18
0.14 18
0.16 18
0.18 18
0.2 18
0.22 18
0.24 18
0.26 18
0.28 18
0.3 18
0.32 18
0.34 18
0.36 18
0.38 18
0.4 18
0.42 18
0.44 18
0.46 18
0.48 18
0.5 18
0.52 18
0.54 18
0.56 18
0.58 18
0.6 18
0.62 18
0.64 18
0.66 18
0.68 18
0.7 18
0.72 18
0.74 18
0.76 18
0.78 18
0.8 18
0.82 18
0.84 18
0.86 18
0.88 18
0.9 18
0.92 18
0.94 18
0.96 18
0.98 18
1 18
};
\addlegendentry{$r_{\delta}$ - $\epsilon_l=1 \times 10^{-8}$}
\addplot [line width=1.2, myblack, dashed]
table {%
0.02 19
0.04 17
0.06 15
0.08 15
0.1 15
0.12 15
0.14 15
0.16 15
0.18 15
0.2 15
0.22 15
0.24 15
0.26 15
0.28 15
0.3 15
0.32 15
0.34 15
0.36 15
0.38 15
0.4 15
0.42 15
0.44 15
0.46 15
0.48 15
0.5 15
0.52 15
0.54 15
0.56 15
0.58 14
0.6 14
0.62 14
0.64 14
0.66 14
0.68 14
0.7 14
0.72 14
0.74 14
0.76 14
0.78 14
0.8 14
0.82 14
0.84 14
0.86 14
0.88 14
0.9 14
0.92 14
0.94 14
0.96 14
0.98 14
1 14
};
\addlegendentry{$r$ - $\epsilon_l=1 \times 10^{-10}$}
\addplot [line width=1.2, myblue, dashed]
table {%
0.02 33
0.04 26
0.06 27
0.08 27
0.1 28
0.12 30
0.14 30
0.16 30
0.18 30
0.2 30
0.22 30
0.24 30
0.26 30
0.28 30
0.3 30
0.32 30
0.34 30
0.36 30
0.38 30
0.4 29
0.42 29
0.44 29
0.46 29
0.48 28
0.5 28
0.52 28
0.54 27
0.56 28
0.58 27
0.6 26
0.62 27
0.64 26
0.66 25
0.68 26
0.7 25
0.72 25
0.74 25
0.76 25
0.78 24
0.8 24
0.82 24
0.84 25
0.86 24
0.88 23
0.9 24
0.92 24
0.94 24
0.96 23
0.98 23
1 23
};
\addlegendentry{$r_{\delta}$ - $\epsilon_l=1 \times 10^{-10}$}
\end{axis}

\end{tikzpicture}}
     }
     \subfigure[]{
         \centering
         \scalebox{.75}{
\begin{tikzpicture}[scale=1]

\definecolor{myblue}{rgb}{0.00000,0.44700,0.74100}%
\definecolor{myred}{rgb}{0.85000,0.32500,0.09800}%
\definecolor{myblack}{rgb}{0,0,0}

\definecolor{gray}{RGB}{128,128,128}
\definecolor{lightgray204}{RGB}{204,204,204}

\begin{axis}[
legend cell align={left},
legend style={
  fill opacity=0.8,
  draw opacity=1,
  text opacity=1,
  at={(0.75,0.4)},
  anchor=center,
  draw=lightgray204
},
log basis y={10},
tick align=outside,
tick pos=left,
x grid style={gray},
xlabel={\(\displaystyle t\)},
xmajorgrids,
xmin=-0.029, xmax=1.049,
xtick style={color=black},
y grid style={gray},
ylabel={\(\displaystyle \mathcal{E}(t)\)},
ymajorgrids,
ymin=3.92800774314104e-09, ymax=3.94686975330337e-06,
ymode=log,
ytick style={color=black}
]
\addplot [line width=1.2, myblack]
table {%
0.02 5.37812061534932e-09
0.04 5.74713121976293e-09
0.06 6.55960130155744e-09
0.08 6.77315930025766e-09
0.1 1.36382016335292e-08
0.12 9.30910977084244e-08
0.14 2.88266776834603e-06
0.16 2.67288833549815e-06
0.18 2.57428783419226e-06
0.2 2.53868504084521e-06
0.22 2.52591803669995e-06
0.24 2.51998407216247e-06
0.26 2.52299904992377e-06
0.28 2.53088899593501e-06
0.3 2.54288338155083e-06
0.32 2.55848291286895e-06
0.34 2.57542627400258e-06
0.36 2.59383875694395e-06
0.38 2.61137192264259e-06
0.4 2.62805272974558e-06
0.42 2.64140474994432e-06
0.44 2.65501668802671e-06
0.46 2.6662427382198e-06
0.48 2.67681657013452e-06
0.5 2.69326548839408e-06
0.52 2.70436383717119e-06
0.54 2.71365059200563e-06
0.56 2.72080062346638e-06
0.58 2.72474248267225e-06
0.6 2.72963455162505e-06
0.62 2.73127441953005e-06
0.64 2.73128978170894e-06
0.66 2.73020261235615e-06
0.68 2.72588074744367e-06
0.7 2.71659147278335e-06
0.72 2.70652264680776e-06
0.74 2.69631021554154e-06
0.76 2.68568494075606e-06
0.78 2.67364213139142e-06
0.8 2.65737925546681e-06
0.82 2.64212411204462e-06
0.84 2.6232014618754e-06
0.86 2.60513816145226e-06
0.88 2.58657866642261e-06
0.9 2.56687379013102e-06
0.92 2.54706538979487e-06
0.94 2.5277018706889e-06
0.96 2.50659546801733e-06
0.98 2.48669353472359e-06
1 2.46706003539136e-06
};
\addlegendentry{$\epsilon_l=1 \times 10^{-6}$}
\addplot [line width=1.2, myblue, dash pattern=on 1pt off 3pt on 3pt off 3pt]
table {%
0.02 5.37812061534932e-09
0.04 5.74713121976293e-09
0.06 6.55960130155744e-09
0.08 6.77315930025766e-09
0.1 1.36382016335292e-08
0.12 9.30910977084244e-08
0.14 1.29107696925594e-07
0.16 1.68078257944287e-07
0.18 2.02461838192774e-07
0.2 2.32823926108658e-07
0.22 2.60718715987861e-07
0.24 2.81778317353944e-07
0.26 3.00441952029147e-07
0.28 3.17272087625824e-07
0.3 3.30267727895155e-07
0.32 3.42465495737609e-07
0.34 3.52592929086258e-07
0.36 3.58096122579025e-07
0.38 3.62579215959212e-07
0.4 3.24838604262904e-07
0.42 3.04892757466556e-07
0.44 3.15220781306918e-07
0.46 3.21666840343788e-07
0.48 3.24954367853419e-07
0.5 3.24966049278982e-07
0.52 3.24501797656815e-07
0.54 3.21166909032907e-07
0.56 3.177363625898e-07
0.58 3.13755791465166e-07
0.6 3.1451168475429e-07
0.62 3.13057253200728e-07
0.64 3.08534420838269e-07
0.66 3.04138728772669e-07
0.68 2.9956173924386e-07
0.7 2.9550465748171e-07
0.72 2.92136502052149e-07
0.74 2.86838233843867e-07
0.76 2.81778771537621e-07
0.78 2.77095556526532e-07
0.8 2.73034223901292e-07
0.82 2.68481363008213e-07
0.84 2.6489498950119e-07
0.86 2.61703394532132e-07
0.88 2.59510528161174e-07
0.9 2.5831354309769e-07
0.92 2.55811336541329e-07
0.94 2.54129637456409e-07
0.96 2.54043020752728e-07
0.98 2.52171395279678e-07
1 2.50907087982156e-07
};
\addlegendentry{$\epsilon_l=1 \times 10^{-8}$}
\addplot [line width=1.2, myred, dashed]
table {%
0.02 5.37812061534932e-09
0.04 5.74713121976293e-09
0.06 6.55960130155744e-09
0.08 7.70084739084093e-09
0.1 7.58214733587928e-09
0.12 8.19666997878906e-09
0.14 8.77801320046741e-09
0.16 8.02881728135665e-09
0.18 7.64275186595382e-09
0.2 6.7724055429939e-09
0.22 6.99636352007189e-09
0.24 6.95549123868225e-09
0.26 6.93265439480107e-09
0.28 6.80929794337479e-09
0.3 6.67827847665985e-09
0.32 6.69111706923406e-09
0.34 6.61741583085291e-09
0.36 6.65148693555112e-09
0.38 6.66625916589984e-09
0.4 6.65654351697353e-09
0.42 6.66334241695694e-09
0.44 6.72751277613835e-09
0.46 6.70416454333703e-09
0.48 6.7049253503272e-09
0.5 6.69239013325757e-09
0.52 6.7183950334005e-09
0.54 6.72109822185285e-09
0.56 6.71933572353677e-09
0.58 6.74553194291255e-09
0.6 6.75931859732466e-09
0.62 6.78358675107998e-09
0.64 6.79735068366322e-09
0.66 6.8096161308734e-09
0.68 6.79489880525457e-09
0.7 6.78572110340841e-09
0.72 6.78722788498703e-09
0.74 6.78169298293846e-09
0.76 6.76721368920808e-09
0.78 6.76067462217519e-09
0.8 6.74973293440047e-09
0.82 6.74257612763963e-09
0.84 6.73114081473626e-09
0.86 6.71642515198712e-09
0.88 6.70076751462537e-09
0.9 6.67835211448539e-09
0.92 6.65774156151556e-09
0.94 6.64367958293674e-09
0.96 6.61970836696125e-09
0.98 6.59936394878216e-09
1 6.57995509269807e-09
};
\addlegendentry{$\epsilon_l=1 \times 10^{-10}$}
\end{axis}

\end{tikzpicture}}
    }
    \caption{Burgers' Equation: Left panel shows the adapted rank over time for various $\epsilon_l$ thresholds. Right panel shows the relative error between TDB-CUR and FOM for different $\epsilon_l$.}
        \label{fig:Burger_rE}
\end{figure}
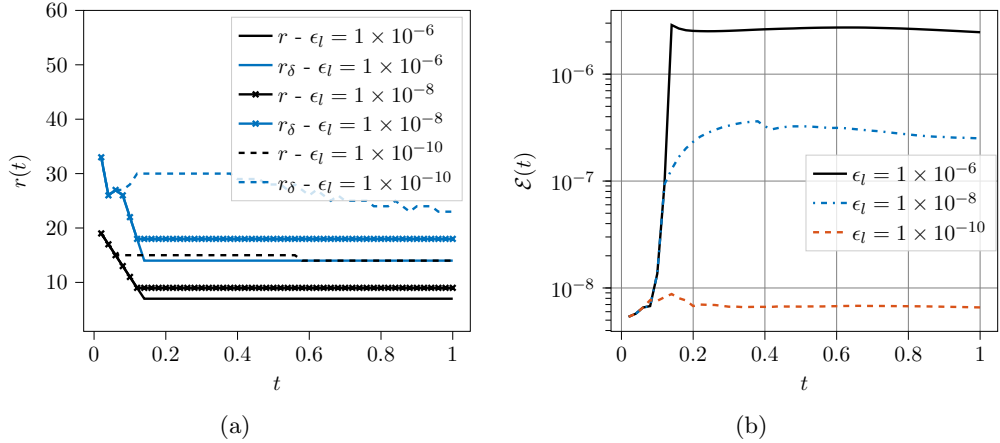

\cref{fig:Burger_rE} shows the effect of varying the rank threshold $\epsilon_l$ and $\epsilon_u$, which controls the accuracy of the TDB-CUR approximation. The left panel shows rank of $\hat{\bm V}$, i.e., $r$, and the rank of Newton's correction matrix, i.e., $ r_{\delta}$, versus time for $\epsilon_l = 10^{-6}, 10^{-8}, 10^{-10}$. For each case, $\epsilon_u$ is set to $10\epsilon_l$, i.e. $\epsilon_u = 10^{-5}, 10^{-7}, 10^{-9}$, respectively. It is clear that $r_{\delta} > r$. Decreasing $\epsilon_l$ results in higher ranks by tightening the accuracy tolerance. Consequently, as shown in the right panel, the relative errors between TDB-CUR and FOM are reduced for smaller values of $\epsilon_l$.
 This demonstrates how the accuracy is controlled by the choice of $\epsilon_l$ and $\epsilon_u$.

\begin{figure}
    \centering
    \scalebox{.57}{
\begin{tikzpicture}[scale=1.2]

\definecolor{myblue}{rgb}{0.00000,0.44700,0.74100}%
\definecolor{myred}{rgb}{0.85000,0.32500,0.09800}%
\definecolor{myblack}{rgb}{0,0,0}%

\definecolor{gray}{RGB}{128,128,128}
\definecolor{lightgray204}{RGB}{204,204,204}

\begin{groupplot}[group style={group size=2 by 1, horizontal sep=2cm}]
\nextgroupplot[
width=4in,
height=3.5in,
legend cell align={left},
legend style={at={(0.98,0.98)},fill opacity=0.8, draw opacity=1, text opacity=1, draw=lightgray204},
log basis y={10},
tick align=outside,
tick pos=left,
title={},
x grid style={gray},
xlabel={Iteration},
xmajorgrids,
xmin= 1, xmax=3.5,
xtick={1,2,3,4},
xtick style={color=black},
y grid style={gray},
ylabel={Residual},
ymajorgrids,
ymin=1e-17, ymax=0.01,
ymode=log,
ytick style={color=black}
]
\addplot [line width=1.2, myblack, mark=x, mark size=3, mark options={solid}]
table {
1 1.85951628e-05
2 1.57217214e-07
3 1.57653946e-07
};
\addlegendentry{AM2 - All entries}
\addplot [line width=1.2, myblack, mark=*, mark size=3, mark options={solid}]
table {
1 1.72659931e-04
2 3.38153299e-08
3 1.27000642e-15
};
\addlegendentry{AM2 - CUR rows}


\addplot [line width=1.2, myblue, mark=x, mark size=3, mark options={solid}]
table {%
1 1.69449732e-05
2 1.28021850e-07
3 1.27927343e-07
};
\addlegendentry{BDF4 - All entries}
\addplot [line width=1.2, myblue, mark=*, mark size=3, mark options={solid}]
table {%
1 1.59644319e-04
2 2.84311347e-08
3 8.27292149e-16
};
\addlegendentry{BDF4 - CUR rows}




\end{groupplot}

\end{tikzpicture}}
    \caption{Burgers' Equation: Average residual versus Newton iterations for different implicit time integration schemes where average residual is equal to $ \| \bm R \|_F/ns$ when is evaluated on all entries and is equal to $\| \bm R(\bm p,:) \|_{F}/ps$  when is evaluated on CUR rows. The results are calculated for $r=5$ and $t=5\Delta t$. }
   \label{fig:Burger_res}
\end{figure}
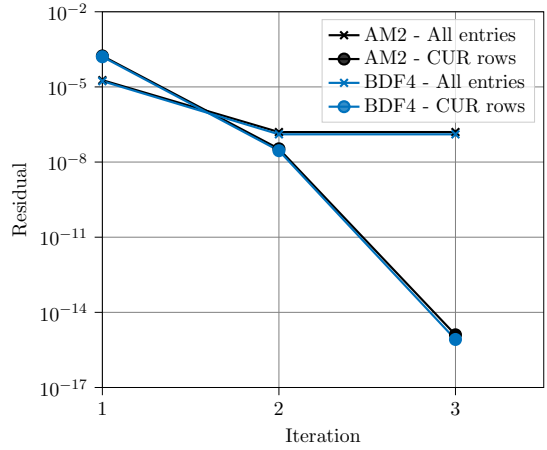
\cref{fig:Burger_res} presents the convergence behavior of the average residual versus the number of Newton iterations for different implicit integrators. It shows the residuals with $e = 0$ for the AM2 and BDF4 methods. The average residual is equal to $\|\bm{R}\|_F/ns$ when evaluated on all entries and $\|\bm{R}(\bm{p},:)\|_F/ps$ when evaluated on CUR DEIM-selected rows. Since no oversampling is used,  the residuals at CUR rows converge quadratically to machine precision.

\begin{figure}
     \centering
     \subfigure[]{
         \centering
\begin{tikzpicture}[scale=0.74]

\definecolor{gray}{RGB}{128,128,128}
\definecolor{lightgray204}{RGB}{204,204,204}
\definecolor{myblue}{rgb}{0.00000,0.44700,0.74100}%
\definecolor{myred}{rgb}{0.85000,0.32500,0.09800}%
\definecolor{mygrey}{rgb}{.7 .7 .7}
\definecolor{myyellow}{rgb}{0.93 0.69 0.13}

\begin{axis}[
legend cell align={left},
legend style={
  fill opacity=0.8,
  draw opacity=1,
  text opacity=1,
  at={(0.05,0.95)},
  anchor=north west,
  draw=lightgray204
},
tick align=outside,
tick pos=left,
x grid style={gray},
xlabel={\(\displaystyle n\)},
xmajorgrids,
xmin=-70.4, xmax=4294.4,
xtick style={color=black},
y grid style={gray},
ylabel={Avg Wall Clock [s]},
ymajorgrids,
ymin=-105.45, ymax=3006.45,
ytick style={color=black}
]
\addplot [line width=1.2, myblue, mark=x, mark size=3.5, mark options={solid}]
table {%
128 43
256 52
512 95
1024 224
2048 695
4096 2865
};
\addlegendentry{FOM}
\addplot [line width=1.2, myred, mark=*, mark size=3.5, mark options={solid}]
table {%
128 36
256 38
512 42
1024 45
2048 53
4096 108
};
\addlegendentry{TDB-CUR}
\end{axis}

\end{tikzpicture}
         \label{fig:Burger_RTN}
     }
     \subfigure[]{
         \centering         
\begin{tikzpicture}[scale=0.74]

\definecolor{gray}{RGB}{128,128,128}
\definecolor{lightgray204}{RGB}{204,204,204}
\definecolor{myblue}{rgb}{0.00000,0.44700,0.74100}%
\definecolor{myred}{rgb}{0.85000,0.32500,0.09800}%
\definecolor{mygrey}{rgb}{.7 .7 .7}
\definecolor{myyellow}{rgb}{0.93 0.69 0.13}

\begin{axis}[
legend cell align={left},
legend style={
  fill opacity=0.8,
  draw opacity=1,
  text opacity=1,
  at={(0.05,0.95)},
  anchor=north west,
  draw=lightgray204
},
tick align=outside,
tick pos=left,
x grid style={gray},
xlabel={\(\displaystyle s\)},
xmajorgrids,
xmin=-70.4, xmax=4294.4,
xtick style={color=black},
y grid style={gray},
ylabel={Avg Wall Clock [s]},
ymajorgrids,
ymin=0, ymax=400.45,
ytick style={color=black}
]
\addplot [line width=1.2, myblue, mark=x, mark size=3.5, mark options={solid}]
table {%
128 9
256 19
512 38
1024 78
2048 157
4096 314
};
\addlegendentry{FOM}
\addplot [line width=1.2, myred, mark=*, mark size=3.5, mark options={solid}]
table {%
128 2
256 4
512 8
1024 15
2048 30
4096 62
};
\addlegendentry{TDB-CUR}
\end{axis}

\end{tikzpicture}
    \label{fig:Burger_RTS}
    }
        \caption{Burgers' Equation: CPU time comparison between FOM and TDB-CUR using AM2, (a) as number of spatial grid points increases, (b) as number of samples increases.}
        \label{fig:Burger_RT}
\end{figure}
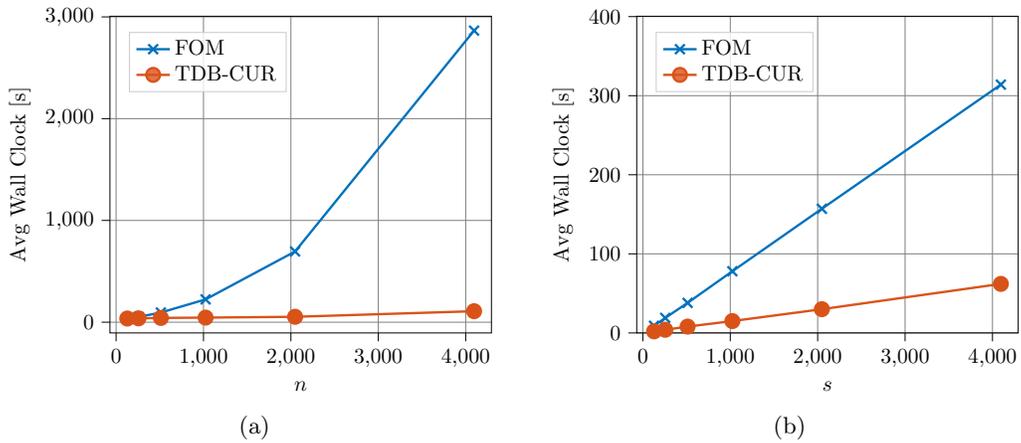

\cref{fig:Burger_RTN} compares the computational cost of FOM and TDB-CUR using AM2 as the number of spatial grid points ($n$) increases, with a fixed number of samples ($s = 512$) and rank ($r = 5$). Here, we use the GMRES method as the linear solver for both FOM and the TDB-CUR methods. The CPU time for FOM increases more than linearly, as the entire $n \times n$ system must be solved. In contrast, TDB-CUR demonstrates linear scaling. Additionally, \cref{fig:Burger_RTS} compares the computational cost of FOM and TDB-CUR using AM2 as the number of samples ($s$) increases, with a fixed number of grid points ($n = 512$) and rank ($r = 5$). The CPU time increases linearly for both FOM and TDB-CUR, with TDB-CUR being significantly more efficient than FOM.

\begin{table}
    \centering
    \begin{tabular}{cccc}
    \hline
         & FOM (s) & TDB-CUR with $r=5$ (s) & TDB-CUR with $r=15$ (s) \\
         \hline
        AM2 & 12.2 & 0.6 & 1.2\\
        IRK2 & 17.4 & 0.7 & 1.4 \\
        IRK4 & 34.2 & 1.8 & 3 \\
        BDF2 & 9.1 & 0.3 & 0.6 \\
         BDF4 & 9.3 & 0.4 & 0.7 \\
         \hline
    \end{tabular}
    \caption{Burgers' Equation: Average computational cost per iteration for the full-order model (FOM) and the TDB-CUR method with ranks $r=5$ and $r=15$, using a setup with $n=2048$ spatial grid points, $s=2048$ samples, and time-step size $\Delta t=0.01$.}
    \label{tab:comp_cost}
\end{table}

To further quantify the computational efficiency of the proposed TDB-CUR method, we analyze the average cost per iteration for different implicit time integration schemes. \cref{tab:comp_cost} presents the average cost of one iteration for the FOM and the TDB-CUR method with ranks $r=5$ and $r=15$, using a setup with $n=2048$ spatial grid points, $s=2048$ samples, and a time-step size of $\Delta t=0.01$. Here, we utilize the GMRES method as the linear solver for both the FOM and computing $\bm V^k(:,\bm s)$ in the TDB-CUR. Across all the examined methods, including Adams-Moulton (AM2), implicit Runge-Kutta (DIRK2 and DIRK4), and backward differentiation formulas (BDF2 and BDF4), the TDB-CUR method exhibits substantial computational speedups compared to the FOM. For instance, with the AM2 method, the TDB-CUR approach is approximately 20 times faster than the FOM for $r=5$, and even with $r=15$, it maintains a significant speedup. Similar trends are observed for the other methods, with the TDB-CUR method being up to 30 times faster than the FOM for $r=5$, while still offering considerable speedups for $r=15$. These results highlight the computational efficiency and scalability of the proposed method, enabling accurate low-rank approximations of nonlinear dynamical systems at a fraction of the computational cost compared to standard full-order simulations.

\subsection{Stochastic 2D Gray-Scott Equations}
For the third example, we consider a Gray-Scott reaction-diffusion equation with a
random coefficient, where we have concentrations of two chemicals $u$, $v$ on a 2D periodic domain. The model is described with the system of PDEs as shown below: 
$$
\begin{aligned}
& \frac{\partial u}{\partial t}=\epsilon_1 (\frac{\partial^2 u}{\partial x^2}+\frac{\partial^2 u}{\partial y^2}) + \alpha(1-u)     - u v^2, && x \in[-1,1], y \in[-1,1], t \in[0,2000]\\
& \frac{\partial v}{\partial t}=\epsilon_2 (\frac{\partial^2 v}{\partial x^2}+\frac{\partial^2 v}{\partial y^2}) - \beta v + u v^2, \\
& u(x,y,0) = 1 - \exp(-80((x + 0.05)^2 + (y + 0.02)^2)), \\
& v(x,y,0) = \exp(-80((x - 0.05)^2 + (y - 0.02)^2)), \\
\end{aligned}
$$
where $\epsilon_1=2\times 10^{-5}$, $\epsilon_2=1\times 10^{-5}$, $\alpha=0.04$. The random coefficient $\beta$ is given by $\beta = 0.1 (1 +  \sigma \xi)$ where $\xi \sim \mathcal{U}(0,1)$ and $\sigma=10^{-4}$. We use a second-order finite difference scheme on a uniform grid with $N_x = 200$ and $N_y = 200$ for discretization of the spatial domain and $s=32$ samples. This example considers the DIRK4 integrator with following parameters: over-sampling $e=25$, timestep $\Delta t=5$, rank thresholds $\epsilon_l = 10^{-12}$, $\epsilon_u = 10^{-11}$, and Newton tolerance $\epsilon_t = 10^{-14}$. 

\begin{figure}
     \centering
     \subfigure[]{
         \centering
         \scalebox{.675}{\includegraphics[]{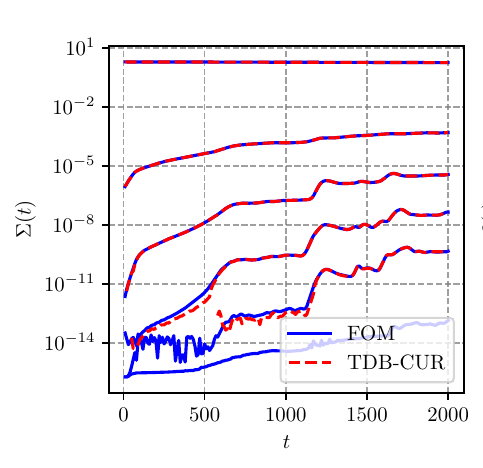}}
         \label{fig:GS-a}
     }
     \subfigure[]{
         \centering
         \scalebox{.675}{\includegraphics[]{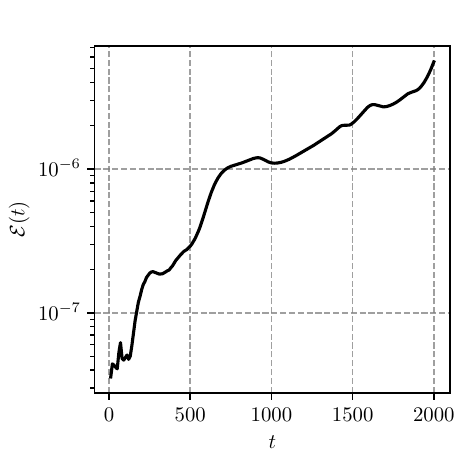}}
         \label{fig:GS-b}
    }
    \vskip\baselineskip
    \subfigure[]{
         \centering
         \scalebox{.7}{\includegraphics[]{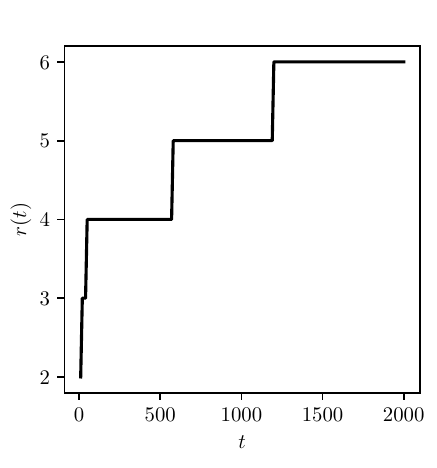}}
         \label{fig:GS-c}
    }
    \subfigure[]{
         \centering
         \scalebox{.7}{\includegraphics[]{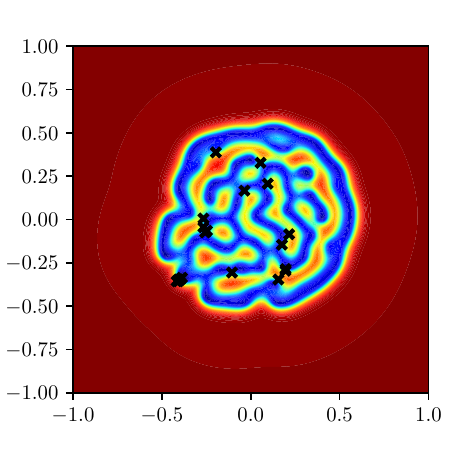}}
         \label{fig:GS-d}
    }
    \caption{2D Gray-Scott Equations: Performance of the implicit TDB-CUR method: (a) Comparison of singular value evolution of FOM and TDB-CUR, (b) Relative error versus time, (c) Adapted rank in different time-steps, (d) Mean solution contour of $u$ at last time-step $t=2000$, where selected DEIM points are shown by black markers.}
        \label{fig:GS}
\end{figure}

\cref{fig:GS-a} displays the evolution of the singular values over time. The leading singular values of the FOM and TDB-CUR match closely, validating the accuracy of the proposed method. \cref{fig:GS-b} shows the relative error between TDB-CUR and FOM over time remains below $10^{-5}$. \cref{fig:GS-c} illustrates the rank increase over time as the dynamics become more complex and require additional modes. Lastly, \cref{fig:GS-d} shows a 2D contour plot of the mean solution for variable $u$ at the final time $t=2000$ across all samples. Also, the DEIM points (selected rows) are shown by black markers.

\begin{figure}[!t]
    \centering
    \includegraphics[width=\textwidth]{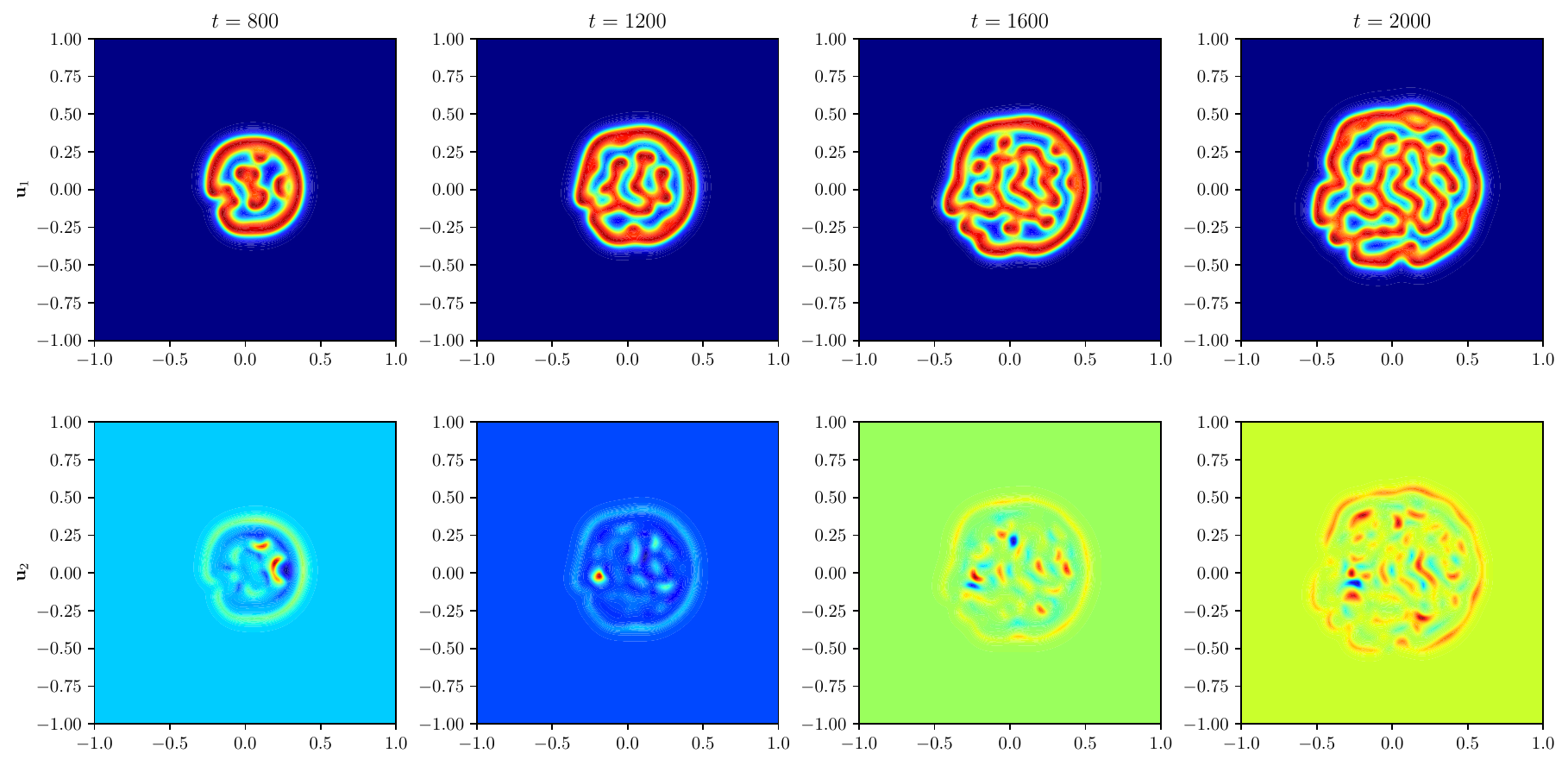}
    \caption{2D Gray-Scott Equations: First two modes of variable $u$ at different time-steps.}
    \label{fig:GS_M}
\end{figure}

\cref{fig:GS_M} shows the evolution of the dominant spatial modes captured by the TDB-CUR model reduction method for the $u$ variable in the stochastic 2D Gray-Scott equations. The figure shows the first two dominant modes associated with the two largest singular values at various time steps throughout the simulation.  As time progresses, the modes become more complex and begin to exhibit the spot/stripe patterns characteristic of the Gray-Scott system. The TDB-CUR method is able to efficiently extract these dominant coherent structures and their evolution over time. Tracking these time-evolving spatial modes is crucial for accurately and efficiently capturing the dynamics with low dimensionality.


\section{\label{sec:Conclusion}Conclusion}
We present a CUR methodology for the implicit time integration of random parametric PDEs on low-rank matrix manifolds. This methodology is computationally efficient because, for a rank-\(r\) approximation, only \(\mathcal{O}(r)\) columns and rows need to be solved. Another advantage of this approach is its applicability to arbitrarily nonlinear PDEs, including those with non-polynomial nonlinearities. The algorithm has been extended to standard high-order implicit time integration methods, including multistep and diagonally implicit Runge-Kutta schemes.

The CUR methodology addresses the implicit nonlinear matrix differential equations (MDEs) at strategically selected columns and rows, which are chosen using the DEIM algorithm or its variants. Newton's method is then employed to solve the nonlinear equations for these selected columns and rows. The implicit time advancement for the columns amounts to independent nonlinear solutions of a deterministic solver for a specific choice of parameters, which can be performed in a non-intrusive manner and parallel. Solving for the selected rows is more complex due to the dependencies among the rows of the MDE resulting from the spatial discretization of differential operators. Therefore, solving only a sparse set of rows is challenging.

To address this, an efficient algorithm has been developed to resolve the row dependencies using a low-rank approximation of the Newton correction matrix. This innovation is crucial for achieving computational efficiency. 

The methodology is demonstrated on analytical problems and PDEs, including stochastic Burgers' and Gray-Scott equations. The results validated the accuracy, robustness, and computational benefits of implicit TDB-CUR compared to full-order model time integration.

\section*{\label{sec:Acknowledgement}Acknowledgement}
 This work is sponsored by a funding from Transformational Tools and Technology (TTT), NASA grant no. 80NSSC22M0282, USA  and by the Air Force Office of Scientific Research award no. FA9550-22-1-0064.

\clearpage
\appendix

\section{\label{Stable_CUR}Stable CUR Algorithm}

The stable CUR pseudocode is presented via \cref{alg:SCUR}  and we refer to \cite{DPNFB23} for more details.

\begin{algorithm}[H]
\caption{Stable CUR Algorithm\label{alg:SCUR}}
\hspace*{\algorithmicindent} \textbf{Input}: $\mathbf{V}(:, \mathbf{s}) \in \mathbb{R}^{n \times r}$, $\mathbf{V}(\mathbf{p},:) \in \mathbb{R}^{r' \times s}$ \\
\hspace*{\algorithmicindent} \textbf{Output}: $\mathbf{U}$, $\boldsymbol \Sigma$, $\mathbf{Y}$
\begin{algorithmic}[1]
\State $\bm{Q},\bm{R} = \texttt{QR}(\mathbf{V}(:, \mathbf{s}))$\Comment{Compute the QR of $\mathbf{V}(:, \mathbf{s})$}
\State $\mathbf{Z}=\mathbf{Q}(\mathbf{p},:)^{\dagger} \mathbf{V}(\mathbf{p},:)$\Comment{Compute $\bm Z$ as an oblique projection of $\bm V$ onto $\bm Q$} 
\State $\mathbf{U}_{\mathbf{Z}}, \boldsymbol \Sigma, \mathbf{Y} = \texttt{SVD}(\mathbf{Z}$)\Comment{Compute the SVD of $\bm Z$}
\State $\mathbf{U}=\mathbf{Q} \mathbf{U}_{\mathbf{Z}}$ \Comment{In-subspace rotation of the orthonormal basis $\bm Q$}
\end{algorithmic}
\end{algorithm}

\bibliographystyle{abbrvnat}
\bibliography{main}

\end{document}